\documentclass[12pt,leqno]{article}
\usepackage{amsmath,amssymb,amscd,latexsym}

\usepackage[all]{xy}
\newcommand{\A}{{\mathbb{A}}}
\newcommand{\C}{{\mathbb{C}}}

\newcommand{\Ge}{\mathbb{G}}
\newcommand{\Pa}{{\mathbb{P}}}
\newcommand{\Q}{{\mathbb{Q}}}
\newcommand{\oQ}{\overline{\Q}}
\newcommand{\hoQ}{\hat{\oQ}}
\newcommand{\Z}{{\mathbb{Z}}}
\newcommand{\oZ}{\overline{\Z}}

\newcommand{\tL}{\tilde{L}}
\newcommand{\Abb}{\mathrm{Alb}}
\newcommand{\abb}{\mathrm{ab}}

\newcommand{\et}{\mathrm{\acute{e}t}}
\newcommand{\ev}{\mathrm{ev}}
\newcommand{\of}{\overline{f}}
\newcommand{\og}{\overline{g}}

\newcommand{\id}{\mathrm{id}}
\renewcommand{\mod}{\;\mathrm{mod}\;}
\newcommand{\rank}{\mathrm{rank}}
\newcommand{\res}{\mathrm{res}}
\newcommand{\spec}{\mathrm{spec}\,}

\newcommand{\Aut}{\mathrm{Aut}}

\newcommand{\Ext}{\mathrm{Ext}}
\newcommand{\uExt}{\underline{\Ext}}

\newcommand{\Gal}{\mathrm{Gal}}
\newcommand{\GL}{\mathrm{GL}\,}
\newcommand{\Hom}{\mathrm{Hom}}
\newcommand{\uHom}{\underline{\Hom}}

\newcommand{\uIso}{\underline{\mathrm{Iso}}}

\newcommand{\Ker}{\mathrm{Ker}\,}
\newcommand{\Lie}{\mathrm{Lie}\,}
\newcommand{\Pic}{\mathrm{Pic}}

\newcommand{\Rep}{\mathbf{Rep}\,}

\newcommand{\Zar}{\mathrm{Zar}}
\newcommand{\tors}{\mathrm{tors}}
\newcommand{\Ah}{{\mathcal A}}

\newcommand{\hAh}{\hat{\Ah}}
\newcommand{\Ch}{{\mathcal C}}
\newcommand{\Fh}{{\mathcal F}}

\newcommand{\Ih}{{\mathcal I}}

\newcommand{\Oh}{{\mathcal O}}

\newcommand{\Yh}{\mathcal{Y}}

\newcommand{\ea}{\mathfrak{a}}

\newcommand{\eo}{\mathfrak{o}}

\newcommand{\eB}{\mathfrak{B}}

\newcommand{\eX}{{\mathfrak X}}

\newcommand{\oK}{\overline{K}}
\newcommand{\ox}{\overline{x}}

\newcommand{\hA}{\hat{A}}
\newcommand{\ha}{\hat{a}}

\newcommand{\hsigma}{\,\!^{\sigma}\!}

\newcommand{\ohne}{\setminus}
\newcommand{\silo}{\stackrel{\sim}{\longrightarrow}}
\newcommand{\tei}{\, | \,}

\newcommand{\verk}{\raisebox{0.03cm}{\mbox{\scriptsize $\,\circ\,$}}}
\newtheorem{theorem}{Theorem}
\newtheorem{lemma}[theorem]{Lemma}
\newtheorem{prop}[theorem]{Proposition}
\newtheorem{defn}[theorem]{Definition}
\newtheorem{cor}[theorem]{Corollary}

\newenvironment{rem}{\noindent {\bf Remark}}{}

\newcommand{\qed}{\mbox{}\hspace*{\fill}$\Box$}

\newenvironment{proof}{\noindent {\bf Proof}}{}
\begin{document}
\title{Line bundles and $p$-adic characters}
\author{Christopher Deninger \\ \small Mathematisches Institut \\[-0.2cm]
\small Einsteinstr. 62, 48149 M\"unster, Germany \\[-0.2cm]
\small deninger@math.uni-muenster.de \and Annette Werner \\ \small Fachbereich Mathematik \\[-0.2cm]
\small Universit\"at Siegen \\[-0.2cm] 
\small Walter-Flex-Str. 3, 57068 Siegen, Germany \\[-0.2cm]
\small email: werner@math.uni-siegen.de}
\date{}
\maketitle
\section{Introduction}
In the paper \cite{De-We2} we defined isomorphisms of parallel transport along \'etale paths for a certain class of vector bundles on $p$-adic curves. In particular, these vector bundles give rise to representations of the fundamental group.

One aim of the present paper is to discuss in more detail the special case of line bundles of degree zero on a curve $X$ with good reduction over $\oQ_p$. By \cite{De-We2} we have a continuous, Galois-equivariant homomorphism
\begin{equation}
  \label{eq:intro1}
  \alpha : \Pic^0_{X / \oQ_p} (\C_p) \longrightarrow \Hom_c (\pi^{\abb}_1 (X) , \C^*_p) \; .
\end{equation}
Here $\Hom_c (\pi^{\abb}_1 (X) , \C^*_p)$ is the topological group of continuous $\C^*_p$-valued characters of the algebraic fundamental group $\pi_1 (X,x)$.

The map $\alpha$ can be rephrased in terms of the Albanese variety $A$ of $X$ as a continuous, Galois-equivariant homomorphism
\begin{equation}
  \label{eq:intro2}
  \alpha : \hA (\C_p) \longrightarrow \Hom_c (TA , \C^*_p) \; .
\end{equation}
Therefore we focus in this paper on abelian varieties $A$ over $\oQ_p$ with good reduction. 

We also consider vector bundles of higher rank on $A_{\C_p} = A \otimes_{\oQ_p} \C_p$. In section 2 we define a category $\eB_{A_{\C_p}}$ of vector bundles on $A_{\C_p}$ which contains all line bundles algebraically equivalent to zero. Then we define for each bundle $E$ in $\eB_{A_{\C_p}}$ a continuous representation $\rho_E$ of the Tate module $TA$ on the fibre of $E$ in the zero section. The association $E \mapsto \rho_E$ is functorial, Galois-equivariant and compatible with several natural operations on vector bundles.

Besides, we show that it is compatible with the theory for curves in \cite{De-We2} via the Albanese morphism $X \to A$. 

Every rank two vector bundle on $A_{\C_p}$ which is an extension of the trivial line bundle $\Oh$ by itself lies in $\eB_{A_{\C_p}}$. The functor $\rho$ induces a homomorphism between $\Ext_{\eB_{A_{\C_p}}} (\Oh , \Oh) \simeq H^1 (A , \Oh) \otimes \C_p$ and the group of continuous extensions
\[
\Ext^1_{TA} (\C_p , \C_p) \simeq H^1_{\et} (X , \Q_p) \otimes \C_p \; .
\]
Hence $\rho$ induces a homomorphism 
\begin{equation}
  \label{eq:intro3}
  H^1 (A , \Oh) \otimes \C_p \longrightarrow H^1_{\et} (X , \Q_p) \otimes \C_p \; .
\end{equation}
In section 3 we show that this homomorphism is the Hodge--Tate map coming from the Hodge--Tate decomposition of $H^1_{\et} (X , \Q_p) \otimes \C_p$. Here we use an explicit description of the Hodge--Tate map by Faltings and Coleman via the universal vectorial extension.

In section 4 we consider the case of line bundles algebraically equivalent to zero on $A$. We prove an alternative description of the homomorphism $\alpha$ in (\ref{eq:intro2}), which shows that the restriction of $\alpha$ to the points of the $p$-divisible group of $\hA$ coincides with a homomorphism defined by Tate in \cite{Ta} using the duality of the $p$-divisible groups associated to $A$ and $\hA$.

In fact, the whole project started with the search for an alternative description of Tate's homomorphism for line bundles on curves which could be generalized to higher rank bundles. 

We prove that $\alpha$ is a $p$-adic analytic morphism of Lie groups whose Lie algebra map
\[
\Lie \alpha : H^1 (A , \Oh) \otimes \C_p = \Lie \hA (\C_p) \longrightarrow \Lie \Hom_c (TA , \C^*_p) = H^1_{\et} (A ,\Q_p) \otimes \C_p
\]
coincides with the Hodge--Tate map (\ref{eq:intro3}).

On the torsion subgroups, $\alpha$ is an isomorphism, so that we get the following commutative diagram with exact lines:
\[
 \vcenter{\xymatrix@-1.2pc{
0 \ar[r] & \hA (\C_p)_{\tors} \ar[r] \ar[d]^{\wr} & \hA (\C_p) \ar[r]^{\log} \ar[d]^{\alpha} & \Lie \hA (\C_p) \ar@{=}[r] \ar[d]^{\Lie \alpha} & H^1 (A , \Oh) \otimes_{\oQ_p} \C_p \ar[r] & 0 \\
0 \ar[r] & \Hom_c (TA , \mu) \ar[r] & \Hom_c (TA , \C^*_p) \ar[r] & \Hom_c (TA , \C_p) \ar@{=}[r] & H^1_{\et} (A , \Q_p) \otimes \C_p \ar[r] & 0
}}
\]
Here $\mu$ is the subgroup of roots of unity in $\C^*_p$. If $A$ is defined over $K$, the vertical maps are all $G_K = \Gal (\oQ_p / K)$-equivariant. 

Besides, we determine the image of $\alpha$. By $CH^{\infty} (TA)$ we denote the group of continuous characters $\chi : TA \to \C^*_p$ whose stabilizer in $G_K$ is open. Then $\alpha$ induces an isomorphism of topological groups between $\hA (\C_p)$ and the closure of $CH^{\infty} (TA)$ in $\Hom_c (TA , \C^*_p)$ with respect to the topology of uniform convergence.

The final section deals with a smooth and proper variety $X$ over $\oQ_p$ whose $H^1$ has good reduction. Combining the previous results for the Albanese variety of $X$ with Kummer theory, we obtain an injective homomorphism of $p$-adic Lie groups
\[
\alpha^{\tau} : \Pic^{\tau}_{X / \oQ_p} (\C_p) \longrightarrow \Hom_c (\pi^{\abb}_1 (X) , \C^*_p)
\]
and determine its image and Lie algebra map.

There is an overlap between parts of the present paper and results independently obtained by Faltings (see \cite{Fa2n}). In \cite{Fa2n}, Faltings develops a more general theory where he proves an equivalence of categories between vector bundles on a curve $X$ over $\C_p$ endowed with a $p$-adic Higgs field and a certain category of ``generalized representations'' of the fundamental group of $X$.

The results of the present paper originally formed the first part of the preprint \cite{De-We1}. However, this preprint will not be published, since the results in its second part are contained in the much more general theory of \cite{De-We2}. 

{\bf Acknowledgements} It is a pleasure to thank Damian Roessler and Peter Schneider for interesting discussions and suggestions.
\section{Vector bundles giving rise to $p$-adic representations}
\label{sec:2}

Let $A$ be an abelian variety over $\oQ_p$ with good reduction, and let $\Ah$ be an abelian scheme over $\oZ_p$ with generic fibre $A$.

We denote the ring of integers in $\C_p = \hoQ_p$ by $\eo$, and put
\[
\eo_n = \eo / p^n \eo = \oZ_p / p^n \oZ_p \; .
\]
We write $\Ah_{\eo} = \Ah \otimes_{\oZ_p} \eo$ and $A_{\C_p} = A \otimes_{\oQ_p} \C_p$. Besides, we denote by $\Ah_n = \Ah \otimes_{\oZ_p} \eo_n$ the reduction of $\Ah$ modulo $p^n$.

Let $\hAh = \Pic^0_{\Ah / \oZ_p}$ be the dual abelian scheme. Its generic fibre $\hA = \Pic^0_{A / \oQ_p}$ is the dual abelian variety of $A$. 

\begin{defn}
  \label{t1}
Let $\eB_{\Ah_{\eo}}$ be the full subcategory of the category of vector bundles on the abelian scheme $\Ah_{\eo}$ consisting of all bundles $E$ on $\Ah_{\eo}$ satisfying the following property:\\
For all $n \ge 1$ there exists some $N = N (n) \ge 1$ such that the reduction $(N^* E)_n$ of $N^* E$ modulo $p^n$ is trivial on $\Ah_n = \Ah \otimes_{\oZ_p} \eo_n$. Here $N : \Ah_{\eo} \to \Ah_{\eo}$ denotes multiplication by $N$.
\end{defn}

Note that every vector bundle $F$ on $A_{\C_p}$ can be extended to a vector bundle $E$ on $\Ah_{\eo}$. This can be shown by induction on the rank of $F$. If $F$ is a line bundle on $A_{\C_p}$, it corresponds to a $\C_p$-valued point in the Picard scheme $\Pic_{A / \oQ_p}$ of $A$. Since every connected component of $\Pic_{A / \oQ_p}$ contains a $\oQ_p$-valued point, there exists a line bundle $M$ on $A$ such that $F \otimes M^{-1}_{\C_p}$ lies in $\Pic^0_{A / \oQ_p} (\C_p)$. As $\Pic^0_{A / \oQ_p} (\C_p) = \Pic^0_{\Ah / \oZ_p} (\eo)$, the line bundle $F \otimes M^{-1}_{\C_p}$ can be extended to $\Ah_{\eo}$. Therefore is suffices to show that every line bundle $M$ on $A$ can be extended to a line bundle on $\Ah$. Now $A$ and $M$ descend to $A_K$ and $M_K$ over a finite extension $K$ of $\Q_p$ in $\oQ_p$ with ring of integers $\eo_K$. The N\'eron model $\Ah_{\eo_K}$ of $A_K$ over the discrete valuation ring $\eo_K$ is noetherian, hence
\[
\Pic_{\Ah_{\eo_K} / \eo_K} (\eo_K) \cong CH^1 (\Ah_{\eo_K}) \longrightarrow CH^1 (A_K) \simeq \Pic_{A_K / K} (K)
\]
is surjective. Therefore $M$ can be extended to $\Ah$.

Now let $F$ be a vector bundle on $A_{\C_p}$. Then there exists a short exact sequence
\[
0 \longrightarrow F_1 \longrightarrow F \longrightarrow F_2 \longrightarrow 0 \; ,
\]
where $F_2$ is a line bundle. By induction, we can assume that $F_1$ and $F_2$ can be extended to vector bundles $E_1$ and $E_2$ on $\Ah_{\eo}$. Flat base change gives an isomorphism
\[
\Ext^1 (E_2 , E_1) \otimes_{\eo} \C_p \silo \Ext^1 (F_2 , F_1) \; ,
\]
hence $F$ can also be extended to $\Ah_{\eo}$.

We are interested in those bundles on $A_{\C_p}$ which have a model in $\eB_{\Ah_{\eo}}$.

\begin{defn}
  \label{t2}
Let $\eB_{A_{\C_p}}$ be the full subcategory of the category of vector bundles on $A_{\C_p}$ consisting of all bundles $F$ on $A_{\C_p}$ which are isomorphic to the generic fibre of a vector bundle $E$ in the category $\eB_{\Ah_{\eo}}$. 
\end{defn}

Consider the Tate modul $TA = \varprojlim A_N (\oQ_p)$, where $A_N (\oQ_p)$ denotes the group of $N$-torsion points in $A (\oQ_p)$.

By $x_{\eo}$ respectively $x_n$ we denote the zero sections on $\Ah_{\eo}$ respectively $\Ah_n$. For a vector bundle $E$ on $\Ah_{\eo}$ we write $E_{x_{\eo}} = x^*_{\eo} E$ viewed as a free $\eo$-module of rank $r = \rank E$. Similarly, we set $E_{x_n} = x^*_n E$ viewed as a free $\eo_n$-module of rank $r$. Note that
\[
E_{x_{\eo}} = \varprojlim_{n} E_{x_n}
\]
as topological $\eo$-modules, if $E_{x_n}$ is endowed with the discrete topology. 

Assume that $E$ is contained in the category $\eB_{\Ah_{\eo}}$ and fix some $n \ge 1$. Then there exists some $N = N (n) \ge 1$, such that the reduction $(N^* E)_n$ is trivial on $\Ah_n$. The structure morphism $\lambda : \Ah_{\eo} \to \spec \eo$ satisfies $\lambda_* \Oh_{\Ah_{\eo}} = \Oh_{\spec \eo}$ universally. Hence $\Gamma (\Ah_n , \Oh) = \eo_n$, and therefore the pullback map
\[
x^*_n : \Gamma (\Ah_n , (N^* E)_n) \silo \Gamma (\spec \eo_n , x^*_n E_n) = E_{x_n}
\]
is an isomorphism of free $\eo_n$-modules. (Note that $N \verk x_n = x_n$.) On $\Gamma (\Ah_n , (N^* E)_n)$ the group $A_N (\oQ_p)$ acts in a natural way by translation. Define a representation $\rho_{E,n} : TA \to \Aut_{\eo_n} (E_{x_n})$ as the composition:
\[
\rho_{E,n} : TA \longrightarrow A_N (\oQ_p) \longrightarrow \Aut_{\eo_n} \Gamma (\Ah_n , (N^* E)_n) \xrightarrow[\mathrm{via}\,x^*_n]{\sim} \Aut_{\eo_n} E_{x_n} \; .
\]

\begin{lemma}
  \label{t3}
For $E$ in $\eB_{\Ah_{\eo}}$ the representations $\rho_{E,n}$ are independent of the choice of $N$ and form a projective system when composed with the natural projection maps $\Aut_{\eo_{n+1}} E_{x_{n+1}} \to \Aut_{\eo_n} E_{x_n}$. 
\end{lemma}

\begin{proof}
  If $N' = N \cdot M$ for some $M \ge 1$, and $(N^* E)_n$ is trivial, it follows from our construction that
\[
\xymatrix{
A_{N'} (\oQ_p) \ar[d]^{\cdot M} \ar[r] & \Aut (E_{x_n}) \\
A_N (\oQ_p) \ar[ur]
}
\]
is commutative. Hence $\rho_{E,n}$ is independent of the choice of $N$.

Fix some $n \ge 1$ and assume that $(N^* E)_{n+1}$ is trivial. Then $(N^* E)_n$ is also trivial, and the natural action of $A_N (\oQ_p)$ on $\Gamma (\Ah_{n+1} , (N^* E)_{n+1})$ induces the natural action of $A_N (\oQ_p)$ on $\Gamma (\Ah_n , (N^* E)_n)$. Hence
\[
\xymatrix{
 & \Aut_{\eo_{n+1}} E_{x_{n+1}} \ar[dd] \\
TA \ar[ur]^{\rho_{E, n+1}} \ar[dr]_{\rho_{E,n}} \\
 & \Aut_{\eo_n } E_{x_n}
}
\]
is commutative, so that $(\rho_{E,n})_{n \ge 1}$ is a projective system. \qed
\end{proof}

By the lemma, we can define for all $E$ in $\eB_{\Ah_{\eo}}$ an $\eo$-linear representation of $TA$ by
\[
\rho_E = \varprojlim_{n} \rho_{E,n} : TA \longrightarrow \Aut_{\eo} (E_{x_{\eo}})\; .
\]
Since each $\rho_{E,n}$ factors over a finite quotient of $TA$, the map $\rho_E$ is continuous, if $\Aut_{\eo} (E_{x_{\eo}}) \simeq \GL_r (\eo)$ for $r = \rank E$ carries the topology induced by the one of $\eo$. 

Note that for any morphism $f : E_1 \to E_2$ of vector bundles in $\eB_{\Ah_{\eo}}$ the natural $\eo_n$-linear map $x^*_n f : (E_1)_{x_n} \to (E_2)_{x_n}$ is $TA$-equivariant with respect to the actions $\rho_{E_1,n}$ and $\rho_{E_2,n}$.
 Hence the association $E \mapsto \rho_E$ defines a functor
\[
\rho : \eB_{\Ah_{\eo}} \longrightarrow \Rep_{TA} (\eo) \; ,
\]
where $\Rep_{TA} (\eo)$ is the category of continuous representations of $TA$ on free $\eo$-modules of finite rank. 

We denote by $x = x_{\eo} \otimes \C_p$ the unit section of $A (\C_p)$.

Let $F$ be a vector bundle in the category $\eB_{A_{\C_p}}$. Then $F$ can be extended to a bundle $E$ in $\eB_{A_{\eo}}$. We define a $\C_p$-linear representation
\[
\rho_F : TA \longrightarrow \Aut_{\C_p} (F_x)
\]
as $\rho_F = \rho_E \otimes_{\eo} \C_p$, where we identify $F_x$ with $E_{x_{\eo}} \otimes_{\eo} \C_p$.

If $E_1$ and $E_2$ are two extensions of $F$ lying in $\eB_{\Ah_{\eo}}$, flat base chance for $H^0$ of the $\Hom$-bundle yields $\Hom_{\Ah_{\eo}} (E_1 , E_2) \otimes_{\eo} \C_p \simeq \Hom_{A_{\C_p}} ((E_1)_{\C_p} , (E_2)_{\C_p})$. 

Let $\varphi : (E_1)_{\C_p} \silo F \silo (E_2)_{\C_p}$ be the identifications of the generic fibres. Then there exists some $m \ge 1$ such that $m \varphi$ is the generic fibre of a morphism $\psi : E_1 \to E_2$. By functoriality, the induced map $x^*_{\eo} \psi : (E_1)_{x_{\eo}} \to (E_2)_{x_{\eo}}$ is $TA$-equivariant. Therefore, identifying $F_x = (E_1)_{x_{\eo}} \otimes_{\eo} \C_p$ and $F_x = (E_2)_{x_{\eo}} \otimes_{\eo} \C_p$, we see that $\rho_{E_1} \otimes_{\eo} \C_p$ and $\rho_{E_2} \otimes_{\eo} \C_p$ coincide. Hence $\rho_F$ is well-defined.

\begin{theorem}
  \label{t4}
i) The category $\eB_{A_{\C_p}}$ is closed under direct sums, tensor products, duals, internal homs and exterior powers. Besides, it is closed under extensions, i.e. if $0 \to F' \to F \to F'' \to 0$ is an exact sequence of vector bundles on $A_{\C_p}$ such that $F'$ and $F''$ are in $\eB_{A_{\C_p}}$, then $F$ is also contained in $\eB_{A_{\C_p}}$. \\
ii) $\eB_{A_{\C_p}}$ contains all line bundles algebraically equivalent to zero. For any bundle in $\eB_{A_{\C_p}}$ the determinant line bundle is algebraically equivalent to zero.\\
iii) The association $F \mapsto \rho_F$ defines an additive exact functor
\[
\rho : \eB_{A_{\C_p}} \longrightarrow \Rep_{TA} (\C_p) \; ,
\]
where $\Rep_{TA} (\C_p)$ is the category of continuous representations of $TA$ on finite-dimensional $\C_p$-vector spaces. This functor commutes with tensor products, duals, internal homs and exterior powers.\\
iv) Let $f : A \to A'$ be a homomorphism of abelian varieties over $\oQ_p$ with good reduction. Then pullback of vector bundles induces an additive exact functor
\[
f^* : \eB_{A'_{\C_p}} \longrightarrow \eB_{A_{\C_p}} \; ,
\]
which commutes with tensor products, duals, internal homs and exterior powers (up to canonical identifications). The following diagram is commutative:
\[
\xymatrix{
\eB_{A'_{\C_p}} \ar[r]^{f^*} \ar[d]_{\rho} & \eB_{A_{\C_p}} \ar[d]^{\rho} \\
\Rep_{TA'} (\C_p) \ar[r]^F & \Rep_{TA} (\C_p)
}
\]
where $F$ is the functor induced by composition with $Tf : TA \to TA'$.
\end{theorem}

\begin{proof}
  i) We only show that $\eB_{A_{\C_p}}$ is closed under extensions, the remaining assertions are straightforward. So consider an extension $0 \to F' \to F \to F'' \to 0$ with $F'$ and $F''$ in $\eB_{A_{\C_p}}$. Fix some $n \ge 1$. Then we find a number $N$ such that $(N^* F')_n$ and $(N^* F'')_n$ are trivial. Hence $(N^* F)_n$ is an extension of two trivial vector bundles. We claim that this implies the triviality of $((p^n N)^* F)_n$. It suffices to show that $(p^n)^*$ induces the zero map on $\Ext^1_{\Ah_n} (\Oh , \Oh) = H^1 (\Ah_n , \Oh)$. The diagram
\[
\xymatrix{
H^1 (\Ah_n , \Oh) \ar[r]^{(p^n)^*} & H^1 (\Ah_n , \Oh)\\
\Lie \Pic^0_{\Ah_n / \eo_n} \ar[r]^{\Lie p^n} \ar[u]_{\wr} & \Lie \Pic^0_{\Ah_n / \eo_n} \ar[u]_{\wr}
}
\]
is commutative by \cite{BLR}, 8.4, Theorem 1. Since $\Lie p^n$ is multiplication by $p^n$ on the $\eo_n$-module $\Lie \Pic^0_{\Ah_n / \eo_n}$, it is zero. Hence $(p^n)^*$ is indeed the zero map on $\Ext^1_{\Ah_n} (\Oh, \Oh)$.\\[0.2cm]
ii) Since $\hAh_{\eo} = \Pic^0_{\Ah_0 / \eo}$ is proper, we have
\[
\Pic^0_{\Ah_{\eo} / \eo} (\eo) = \Pic^0_{A_{\C_p} / \C_p} (\C_p) \; ,
\]
so that any line bundle $L_{\C_p}$ on $A_{\C_p}$ which is algebraically equivalent to zero, can be extended to a line bundle $L$ on $\Ah_{\eo}$ giving rise to a class in $\hAh_{\eo} (\eo)$. \\
$\hAh$ descends to an abelian scheme $\hAh_{\eo_K}$ over the ring of integers in some finite extension $K$ of $\Q_p$ in $\oQ_p$. The ring $\eo_n = \eo / p^n \eo = \oZ_p / p^n \oZ_p$ is the union of the finite rings $\eo_L / p^n \eo_L$, where $L$ runs through the finite extensions of $K$ in $\oQ_p$. Therefore $\hAh_{\eo} (\eo_n) = \hAh_{\eo_K} (\eo_n)$ is the union of all $\hAh_{\eo_K} (\eo_L / p^n \eo_L)$. Since $\hAh_{\eo_K}$ is of finite type over $\eo_K$, all these groups $\hAh_{\eo_K} (\eo_L / p^n \eo_L)$ are finite. Hence $\hAh_{\eo} (\eo_n)$ is a torsion group.\\
In particular, we find some $N$ such that $N$ annihilates the class of $L_n$ in $\hAh_{\eo} (\eo_n)$. Then $(N^* L)_n$ is trivial, which shows that $L$ is contained in $\eB_{A_{\C_p}}$. \\
If $E$ is a vector bundle in $\eB_{\Ah_{\eo}}$, then by i) its determinant line bundle $L$ is also contained in $\eB_{\Ah_{\eo}}$. Hence there exists some $N \ge 1$ such that $N^* L_1$ is trivial on $\Ah_1$, where $L_1$ and $\Ah_1$ denote the reductions modulo $p$. If $k$ denotes the residue field of $\eo$, this implies that $N^* L_k$ is trivial on $\Ah_k$. Since the N\'eron--Severi group of $\Ah_k$ is torsion free, $L_k$ lies in $\Pic^0_{A_k / k} (k)$.\\
Note that $\Ah$ is projective by \cite{Ray}, Th\'eor\`eme XI 1.4. Hence $\Pic^0_{\Ah_{\eo} / \eo}$ is an open subscheme of the Picard scheme $\Pic_{\Ah_{\eo} / \eo}$. Since the reduction of the point in $\Pic$ induced by $L$ lies in $\Pic^0$, the generic fibre of $L$ is also contained in $\Pic^0$, whence our claim.\\[0.2cm]
iii) The fact that $F \mapsto \rho_F$ is functorial on $\eB_{A_{\C_p}}$ follows from the fact that the corresponding association $E \mapsto \rho_E$ is functorial on $\eB_{\Ah_0}$. The remaining claims in iii) are straightforward.\\[0.2cm]
iv) It is clear that $f^*$ induces a functor $f^* : \eB_{A'_{\C_p}} \to \eB_{A_{\C_p}}$ with the claimed properties. In order to show that the desired diagram commutes, it suffices to show that
\[
\xymatrix{
\eB_{\Ah'_{\eo}} \ar[r]^{f^*} \ar[d]_{\rho} & \eB_{\Ah_{\eo}} \ar[d]^{\rho} \\
\Rep_{TA} (\eo) \ar[r]^F & \Rep_{TA} (\eo)
}
\]
commutes, where $f : \Ah_{\eo} \to \Ah'_{\eo}$ comes from the canonical extension of $f$ to the N\'eron models. Let $E'$ be an object in $\eB_{\Ah'_{\eo}}$. If $x'_{\eo}$ is the zero section of $\Ah'_{\eo}$, we have $f (x_{\eo}) = x'_{\eo}$, so that there is a canonical identification of $E'_{x'_{\eo}}$ and $(f^* E')_{x_{\eo}}$. For every $a \in A_N (\oQ_p)$ we have $t_{f (a)} \verk f = f \verk t_a$, where $t_{f (a)}$ and $t_a$ are the translation maps. If we now go through the definition of the representation of $TA$ on $(f^* E')_{x_{\eo}}$, our claim follows. \qed
\end{proof}

Let $K$ be a finite extension of $\Q_p$ in $\oQ_p$ so that $A$ is defined over $K$, i.e. $A = A_K \otimes_K \oQ_p$ for an abelian variety $A_K$ over $K$. Then $G_K = \Gal (\oQ_p / K)$ acts in a natural way on $TA$ and on $\C_p$ and hence on the category $\Rep_{TA} (\C_p)$ by putting both actions together. To be precise, for $\sigma \in G_K$ put $\hsigma V = V \otimes_{\C_{p, \sigma}} \C_p$ for every finite-dimensional $\C_p$-vector space $V$, and write $\sigma : V \to \hsigma V$ for the natural $\sigma$-linear map. For every representation $\varphi : TA \to \Aut_{\C_p} (V)$, we define $\sigma_* \varphi$ as the representation 
\[
\sigma_* \varphi : TA \xrightarrow{\sigma^{-1}} TA \xrightarrow{\varphi} \Aut_{\C_p} V \xrightarrow{c_{\sigma}} \Aut_{\C_p} \hsigma V \; ,
\]
where $c_{\sigma} (f) = \sigma \verk f \verk \sigma^{-1}$.

For every vector bundle $F$ in $\eB_{A_{\C_p}}$ the vector bundle $\hsigma F = F \times_{\spec \C_p , \spec \sigma} \spec \C_p$ is also contained in the category $\eB_{A_{\C_p}}$. 

\begin{prop}
\label{t5}
For every $\sigma \in G_K$ the following diagram is commutative:
\[
\xymatrix{
\eB_{A_{\C_p}} \ar[r]^-{\rho} \ar[d]^{g_{\sigma}} & \Rep_{TA} (\C_p) \ar[d]^{\sigma_*} \\
\eB_{A_{\C_p}} \ar[r]^-{\rho} & \Rep_{TA} (\C_p)
}
\]
where the functor $g_{\sigma}$ maps $F$ to $\hsigma F$.
\end{prop}

\begin{proof}
  This can be checked directly. \qed
\end{proof}

Let $X$ be a smooth, connected, projective curve over $\oQ_p$. We fix a point $x \in X (\oQ_p)$ and denote by $\pi_1 (X,x)$ the algebraic fundamental group with base point $x$.

In \cite{De-We2}, we define and investigate the category $\eB_{X_{\C_p}}$ of all vector bundles $F$ on $X_{\C_p}$ which can be extended to a vector bundle $E$ on $\eX_{\eo} = \eX \otimes_{\oZ_p} \eo$ for a finitely presented, flat and proper model $\eX$ of $X$ over $\oZ_p$ and which have the following property: For all $n \ge 1$ there exists a finitely presented proper $\oZ_p$-morphism $\pi : \Yh \to \eX$ with finite, \'etale generic fibre such that $\pi^*_n E_n$ is trivial. Here $\pi_n$ and $E_n$ denote again the reductions modulo $p^n$.

Besides, we define in \cite{De-We2} a functor $\rho$ from $\eB_{X_{\C_p}}$ to the category of continuous representations of the \'etale fundamental groupoid of $X$. In particular, every bundle $F$ in $\eB_{X_{\C_p}}$ induces a continuous representation $\rho_F : \pi_1 (X,x) \to \Aut_{\C_p} (F_x)$ of the fundamental group, c.f. \cite{De-We2}, Proposition 20. Let $\Rep_{\pi_1 (X,x)} (\C_p)$ be the category of continuous representations of $\pi_1 (X,x)$ on finite-dimensional $\C_p$-vector spaces. Then the association $F \mapsto \rho_F$ induces a functor $\rho : \eB_{X_{\C_p}} \to \Rep_{\pi_1 (X,x)} (\C_p)$. 

Let us now assume that $X$ has good reduction, i.e. there exists a smooth, proper, finitely presented model $\eX$ of $X$ over $\oZ_p$. Then $\Ah = \Pic^0_{\eX / \oZ_p}$ is an abelian scheme. Let $A = \Pic^0_{X / \oQ_p}$ be the Jacobian of $X$, and let $f : X \to A$ be the embedding mapping $x$ to $0$. After descending to a suitable finite extension of $\Q_p$ in $\oQ_p , \Ah$ becomes the N\'eron model of $A$. Hence $f$ has an extension $f_0 : \eX \to \Ah$. Besides, $f$ induces a homomorphism
\[
f_* : \pi_1 (X,x) \longrightarrow \pi_1 (A, 0) = TA \; ,
\]
which identifies $TA$ with the maximal abelian quotient of $\pi_1 (X,x)$. 

Let $F$ be a vector bundle in $\eB_{A_{\C_p}}$, and let $E$ be a model of $F$ on $\Ah_{\eo}$ such that for all $n \ge 1$ there exists some $N \ge 1$ satisfying $(N^* E)_n$ trivial. Consider the pullback of the covering $N : \Ah \to \Ah$ to $\eX$, i.e. the Cartesian diagram
\[
\xymatrix{
\Yh \ar[r]^g \ar[d]_{\pi_N} & \Ah \ar[d]^N \\
\eX \ar[r]^f & \Ah
}
\]
Then $(\pi^*_N f^* E)_n$ is also trivial, and we see that $f^* F$ is contained in $\eB_{X_{\C_p}}$. 

Hence pullback via $f$ induces a functor $f^* : \eB_{A_{\C_p}} \to \eB_{X_{\C_p}}$.

\begin{lemma}
  \label{t6}
The following diagram is commutative:
\[
\xymatrix{
\eB_{A_{\C_p}} \ar[r]^{f^*} \ar[d]_{\rho} & \eB_{X_{\C_p}} \ar[d]^{\rho} \\
\Rep_{TA} (\C_p) \ar[r]^{\tilde{f}} & \Rep_{\pi_1 (X,x)} (\C_p)
}
\]
where $\tilde{f}$ is the functor induced by composition with the homomorphism $f_* : \pi_1 (X,x) \to TA$.
\end{lemma}

\begin{proof}
  One argues similarly as in the proof of Theorem \ref{t4}, iv). \qed
\end{proof}


\section{The Hodge--Tate map}
\label{sec:3}

In this section we show that the functor $\rho$ can be used to describe the Hodge--Tate decomposition of the first \'etale cohomology of $A$, when we apply it to extensions of the trivial line bundle with itself.

Let $K$ be a finite extension of $\Q_p$ in $\oQ_p$ such that there is an abelian variety $A_K$ over $K$ with $A = A_K \otimes_K \oQ_p$. Put $G_K = \Gal (\oQ_p / K)$.

The Hodge--Tate decomposition (originating from \cite{Ta})
\begin{equation}
  \label{eq:1}
  H^1_{\et} (A, \Q_p) \otimes \C_p \simeq (H^1 (A_K , \Oh) \otimes_K \C_p) \oplus (H^0 (A_K , \Omega^1) \otimes_K \C_p (-1))
\end{equation}
gives rise to a $G_K$-equivariant map
\begin{equation}
  \label{eq:2}
  \theta^*_A : H^1 (A_K , \Oh) \otimes_K \C_p \longrightarrow H^1_{\et} (A , \Q_p) \otimes \C_p \; .
\end{equation}
As Faltings (\cite{Fa2}, Theorem 4) and Coleman (\cite{Co1}, p. 379 and \cite{Co3}, \S\,4) have shown, $\theta^*_A$ has the following elegant description.

Consider the universal vectorial extension of $\Ah$ over the ring $\oZ_p$
\[
0 \longrightarrow \omega_{\hAh} \longrightarrow E \longrightarrow \Ah \longrightarrow 0 \; .
\]
Here $\omega_{\hAh}$ is the vector group induced by the invariant differentials on $\hAh$, i.e. $\omega_{\hAh} (S) = H^0 (S , e^* \Omega^1_{\hAh_S / S})$ for all $\oZ_p$-schemes $S$, where $e$ denotes the zero section.

For $\nu \ge 1$ consider the map
\[
\Ah_{p^{\nu}} (\eo) \longrightarrow \omega_{\hAh} (\eo) / p^{\nu} \omega_{\hAh} (\eo)
\]
obtained by sending $a_{p^{\nu}}$ to $p^{\nu} b_{p^{\nu}}$, where $b_{p^{\nu}} \in E (\eo)$ is a lift of $a_{p^{\nu}}$. Passing to inverse limits, we get a $\Z_p$-linear homomorphism
\[
\theta_A : T_p A \longrightarrow \omega_{\hAh} (\eo) \; .
\]
The $\C_p$-dual of the resulting map
\[
\theta_A : T_p A \otimes \C_p \longrightarrow \omega_{\hAh} (\eo) \otimes \C_p = \omega_{\hAh} (\C_p)
\]
is the map $\theta^*_A$ in (\ref{eq:1}).

In \cite{Co1}, Coleman proved that together with a map defined by Fontaine in \cite{Fo}
\[
H^0 (A , \Omega^1) \otimes \C_p (-1) \longrightarrow H^1_{\et} (A , \Q_p) \otimes \C_p \; ,
\]
$\theta^*_A$ gives the Hodge--Tate decomposition.

Let us write $\Ext^1_{\eB_{A_{\C_p}}} (\Oh , \Oh)$ for the Yoneda group of isomorphism classes of extensions $0 \to \Oh \to \Oh (F) \to \Oh \to 0$, where $F$ is a vector bundle in $\eB_{A_{\C_p}}$ with sheaf of sections $\Oh (F)$ and $\Oh = \Oh_{A_{\C_p}}$. By Theorem \ref{t4}, $\eB_{A_{\C_p}}$ contains all vector bundles which are extensions of the trivial bundle by itself. Hence $\Ext^1_{\eB_{A_{\C_p}}} (\Oh , \Oh)$ coincides with the group $\Ext^1 (\Oh , \Oh)$ in the category of locally free sheaves on $A_{\C_p}$, so that
\[
\Ext^1_{\eB_{A_{\C_p}}} (\Oh , \Oh) = \Ext^1 (\Oh , \Oh) = H^1 (A_{\C_p} , \Oh) = H^1 (A_K , \Oh) \otimes_K \C_p \; .
\]
Since the functor $\rho$ is exact by Theorem \ref{t4}, it induces a homomorphism of $\Ext$ groups
\[
\rho_* : \Ext^1_{\eB_{A_{\C_p}}} (\Oh , \Oh) \longrightarrow \Ext^1_{\Rep_{TA} (\C_p)} (\C_p , \C_p) \; .
\]
There is a natural isomorphism
\[
\Ext^1_{\Rep_{TA} (\C_p)} (\C_p , \C_p) \simeq \Hom_c (TA , \C_p) \; ,
\]
where $\Hom_c (TA , \C_p)$ denotes the continuous homomorphisms from $TA$ to $\C_p$, which is defined as follows. For every extension $0 \to \C_p \xrightarrow{i} V \xrightarrow{\varepsilon} \C_p \to 0$ in $\Rep_{TA} (\C_p)$ choose any $v \in V$ with $\varepsilon (v) = 1$ and define $\psi : TA \to \C_p$ by $\psi (\gamma) = i^{-1} (\gamma v - v)$. Since $\Hom_c (TA, \C_p) = \Hom_{\Z_p} (T_p A , \C_p)$, we get an isomorphism
\[
\Ext^1_{\Rep_{TA} (\C_p)} (\C_p , \C_p) \simeq H^1_{\et} (A , \Q_p) \otimes \C_p \; .
\]

\begin{theorem}
  \label{t7}
The following diagram commutes:
\begin{equation}
  \label{eq:3}
\vcenter{  \xymatrix{
\Ext^1_{\eB_{A_{\C_p}}} (\Oh , \Oh) \ar[r]^-{\rho_*} \ar@{=}[d] & \Ext^1_{\Rep_{TA} (\C_p)} (\C_p , \C_p) \ar[d]^{\wr} \\
H^1 (A_K , \Oh) \otimes_K \C_p \ar[r]^{\theta^*_A} & H^1_{\et} (A , \Q_p) \otimes \C_p
}}
\end{equation}
where $\theta^*_A$ is the map (\ref{eq:2}) appearing in the Hodge--Tate decomposition.
\end{theorem}

\begin{rem}
 Theorem \ref{t7} gives the following novel construction of the Hodge Tate map $\theta^*_A$. Consider a class $c$ in $H^1 (\Ah _{\eo}, \Oh)$. It can be viewed as an extension 
\[
0 \to \Oh \to \Oh (E) \to \Oh \to 0
\]
of locally free sheaves on $\Ah_{\eo}$. The bundle $E$ lies in $\eB_{\Ah_{\eo}}$. Hence, for every $n \ge 1$ there is some $N \ge 1$, in fact, $N = p^n$ will do, such that $(N^* E)_n$ is the trivial rank two bundle on $\Ah_n$. The short exact sequence $0 \to \eo_n \to E_{x_n} \to \eo_n \to 0$ of fibres along the zero section of $\Ah_n$ becomes $TA$-equivariant if $TA$ acts trivially on the $\eo_n$'s and via the projection $TA \to A_N (\oK)$ and the isomorphism
\[
\Gamma (\Ah_n , (N^* E)_n ) \overset{\res}{\silo} E_{x_n}
\]
on $E_{x_n}$. Passing to projective limits, we get a short exact sequence \\
$0 \to \eo \xrightarrow{i} E_x \xrightarrow{q} \eo \to 0$ of $TA$-modules. Set $g_1 = i (1)$ and choose $g_2 \in E_x$ such that $q (g_2) = 1$. Then $E_x$ is a free $\eo$-module on $g_1$ and $g_2$, and the action of $\gamma \in TA$ on $E_x$ is given in terms of the basis $g_1, g_2$ by a matrix of the form $\left( 
  \begin{smallmatrix}
    1 & \beta (\gamma) \\ 0 & 1
  \end{smallmatrix} \right)$ where $ \beta : TA \to \eo$ is a continuous homomorphism. Note that $\beta$ does not depend on the choice of $e_1$. Viewing $\beta$ as an element of $H^1_{\et} (A , \Z_p) \otimes \eo$ we have $\theta^*_A (c) = \beta$.

\end{rem}

\begin{proof}
  Since $H^1 (\Ah_{\eo} , \Oh) \otimes_{\eo} \C_p = H^1 (A_{\C_p}, \Oh)$, we find for every element in \\
$\Ext^1_{\eB_{A_{\C_p}}} (\Oh , \Oh)$ a $p$-power multiple lying in $\Ext^1_{\eB_{\Ah_{\eo}}} (\Oh , \Oh)$. Since $\rho_*$ is a homomorphism between Yoneda $\Ext$-groups, it suffices to show that
\[
\xymatrix{
\Ext^1_{\eB_{\Ah_{\eo}}} (\Oh , \Oh) \ar[r]^{\rho_*} \ar@{=}[d] & \Ext^1_{\Rep_{TA} (\eo)} (\eo , \eo) \ar[d]^{\wr} \\
H^1 (\Ah_{\eo} , \Oh) \ar[r]^{\theta^*_A} & \Hom_c(T_p A , \eo)
}
\]
commutes. Consider an extension $0 \to \Oh \stackrel{i}{\rightarrow} \Fh \to \Oh \to 0$ on $\Ah_\eo$, where $\Fh = \Oh(E)$ for some $E$  in $\eB_{\Ah_{\eo}}$. 

Let $\Ih = \uIso_{\Ext} (\Oh^2 , \Fh)$ be the (Zariski) sheaf on $\Ah_{\eo}$ which associates to an open subset $U \subseteq \Ah_{\eo}$ the set of isomorphisms $\varphi : \Oh^2_U \to \Fh_U$ of extensions, i.e. such that
\[
\xymatrix{
0 \ar[r] & \Oh_U \ar[r] \ar@{=}[d] & \Oh^2_U \ar[r]^p \ar[d]_{\varphi} & \Oh_U \ar[r] \ar@{=}[d] & 0 \\
0 \ar[r] & \Oh_U \ar[r]^i & \Fh_U \ar[r] & \Oh_U \ar[r] & 0
}
\]
commutes. Then $c \in \Ge_a(U) = \Gamma (U, \Oh)$ acts in a natural way on $\Ih (U)$ by mapping
\[
\varphi \longmapsto \varphi + i \verk f_c \verk p \; , 
\]
where $f_c : \Oh_U \to \Oh_U$ is multiplication by $c$. Note that $\Ih$ is a $\Ge_a$-torsor on $\Ah_{\eo}$ and that the class $[\Ih]$ of $\Ih$ in $H^1 (\Ah_{\eo} , \Oh) = H^1_{\Zar} (\Ah_{\eo} , \Ge_a)$ coincides with the image of the extension class given by $\Fh$ under the isomorphism 
\[
\Ext^1 (\Oh , \Oh) \silo H^1 (\Ah_{\eo} , \Oh) \; .
\]
The association
\begin{equation}
  \label{eq:4}
(T \xrightarrow{t} \Ah_{\eo}) \longmapsto \uIso_{\Ext} (\Oh^2_T , t^* \Fh)
\end{equation}
also defines a sheaf on the  flat site over $\Ah_{\eo}$, which is represented by a $\Ge_a$-torsor $Z \to \Ah_{\eo}$. 

We have $H^1_{\Zar} (\Ah_{\eo} , \Ge_a) = H^1_{\rm fppf} (\Ah_{\eo} , \Ge_a) = \Ext^1 (\Ah_{\eo} , \Ge_a)$, so that $Z$ can be endowed with a group structure sitting in an extension of $\Ah_{\eo}$ by $\Ge_a$:
\[
0 \longrightarrow \Ge_a \stackrel{j}{\longrightarrow} Z \longrightarrow \Ah_{\eo} \longrightarrow 0 \; .
\]
Hence there is a homomorphism $h : \omega_{\hAh_{\eo}} \to \Ge_a$ such that $Z$ is the pushout of the universal vectorial extension:
\[
\xymatrix{
0 \ar[r] & \omega_{\hAh_{\eo}} \ar[r] \ar[d]_h & E_{\eo} \ar[r] \ar[d] & \Ah_{\eo} \ar[r] \ar@{=}[d] & 0 \\
0 \ar[r] & \Ge_{a} \ar[r] & Z \ar[r] & \Ah_{\eo} \ar[r] & 0 \; .
}
\]
Recall that $\theta_A : T_p A \to \omega_{\hAh} (\eo)$ is defined as the limit of the maps
\[
\theta_{A,n} : \Ah_{p^n} (\oZ_p) \longrightarrow \omega_{\hAh} (\eo) / p^n \omega_{\hAh} (\eo)
\]
associating to $a \in \Ah_{p^n} (\oZ_p)$ the class of $p^n b$ for an arbitrary preimage $b \in E_{\eo} (\eo)$ of $a$.

By \cite{Ma-Me}, chapter I, the natural isomorphism
\[
\Hom_{\eo}(\omega_{\hAh} (\eo) , \eo) \stackrel{\sim}{\longrightarrow} \mbox{Lie}(\hAh_{\eo})(\eo) \stackrel{\sim}{\longrightarrow}  H^1 (\Ah_{\eo} , \Oh) = \mbox{Ext}^1(\Ah_\eo, \Ge_{a})
\]
sends a map to the corresponding pushout of $E_{\eo}$.  Hence $\theta^*_A : H^1 (\Ah_{\eo} , \Oh) \to \Hom_c (T_p A , \eo)$ maps the extension class of $\Fh$ to the limit of the maps
\[
\Ah_{p^n} (\oZ_p) \xrightarrow{\theta_{A,n}} \omega_{\hAh} (\eo) / p^n \omega_{\hAh} (\eo) \xrightarrow{h_n} \eo_n \; , 
\]
where $h_n$ is induced by $h$.

This map can also be described as follows: For $a \in \Ah_{p^n} (\oZ_p) \subseteq \Ah_{p^n} (\eo)$ choose a preimage $z \in Z (\eo)$. Then
\[
h_n \verk \theta_{A,n} (a) = \mbox{class of} \; p^n z \; \mbox{in} \; \Ge_a(\eo) / p^n \Ge_a(\eo) = \eo_n \; .
\]
Set $Z_n = Z \otimes_{\eo} \eo_n$. Let $\pi_{p^n}$ denote multiplication by $p^n$ on $\Ah_n$. We have seen in the proof of Theorem \ref{t4}, i) that $\pi_{p^n}^* \Fh_n$ is a trivial extension. Hence $\pi_{p^n}^* Z_n$ is trivial in $\Ext(\Ah_n, \Ge_{a, \eo_n})$, and there is a splitting   $r : \Ah_n \to \pi_{p^n}^* Z_n$ of the extension
\[
0 \longrightarrow \Ge_{a , \eo_n} \longrightarrow \pi_{p^n}^* Z_n \longrightarrow \Ah_n \longrightarrow 0
\]
over $\eo_n$. Let $g : \pi_{p^n}^* Z_n \rightarrow Z_n$ denote the projection, and denote by 
$a_n \in \Ah_n(\eo_n)$ the point induced by $a$. Then $g (r(a_n))$ projects to zero in $\Ah_n$, hence it is equal to $j(c)$ for some $c \in \Ge_a(\eo_n) = \eo_n$. Since for any $a' \in
 \Ah_n(\eo_n)$ with $p^n a' = a_n $ the point $g (r(a'))$ is a preimage of $a_n$, we have $h_n \verk \theta_{A,n} (a) = c \in \eo_n$. Besides,  $Z$ represents the functor (\ref{eq:4}), so that the map $r: \Ah_n \rightarrow \pi_{p^n}^* Z_n$  corresponds to a trivialization
\[
\xymatrix{
0 \ar[r] & \Oh_{\Ah_n} \ar[r] \ar@{=}[d] & \Oh^2_{\Ah_n} \ar[r] \ar[d]_{\varphi} & \Oh_{\Ah_n} \ar[r] \ar@{=}[d] & 0 \\
0 \ar[r] & \Oh_{\Ah_n} \ar[r] & \pi_{p^n}^* \Fh_n \ar[r] & \Oh_{\Ah_n} \ar[r] & 0 \; .
}
\]
The point $g (r(a_n))$ in the kernel of $Z_n(\eo_n) \rightarrow \Ah_n(\eo_n)$ corresponds to the
trivialization
\[
\alpha: \Oh^2_{\spec \eo_n} \stackrel{a_n^* \varphi}{\longrightarrow} a_n^* \pi_{p^n}^* \Fh  \stackrel{\sim}{\longrightarrow} 0^* \Fh,
\]
where $0$ is the zero element in $\Ah_n(\eo_n)$. 
Besides, the trivialization 
\[
\beta: \Oh^2_{\spec \eo_n} \stackrel{0^* \varphi}{\longrightarrow} 0^* \pi_{p^n}^* \Fh  \stackrel{\sim}{\longrightarrow} 0^* \Fh
\]
is given by the zero element in $Z_n$. By definition, $\alpha = \beta + i \circ f_c \circ p$.
If we denote the canonical basis of $\Gamma( \Ah_n, \Oh^2_{\Ah_n})$ by $e_1,e_2$, and the induced basis of $\Gamma(\Ah_n, \pi_{p^n}^* \Fh)$ by  $f_1 , f_2$, it follows that
$ a_n^* f_2 - 0^* f_2 = i(c)$. 

On the other hand, the image of $E$ under 
\[
\rho_* : \Ext^1_{\eB_{\Ah_{\eo}}} (\Oh , \Oh) \longrightarrow \Ext^1_{\Rep_{TA} (\eo)} (\eo , \eo) = \Hom_c(T_p A , \eo)\]
is the  homomorphism $\gamma: T_pA \rightarrow \eo$, such that $\gamma \mod p^n$ maps the $p^n$-torsion point
$a_n \in \Ah_n(\eo_n)$ to the element 
\[i^{-1}( 0^*(\tau_{a_n}^* f_2)- 0^* f_2) = i^{-1}(a_n^* f_2 - 0^*f_2),\] 
where $\tau_{a_n}$ is translation by $a_n$, and hence to $c = h_n \circ \theta_{A,n}(a)$.
This proves our claim. \qed
\end{proof}

\begin{cor}
  \label{t8}
Let $X$ be a smooth, connected, projective curve over $\oQ_p$ with good reduction and let $x \in X (\oQ_p)$ be a base point. The functor $\eB_{X_{\C_p}} \to \Rep_{\pi_1 (X,x)} (\C_p)$ from \cite{De-We2} induces a homomorphism
\[
\rho_* : \Ext^1_{\eB_{X_{\C_p}}} (\Oh , \Oh) \longrightarrow \Ext^1_{\Ch} (\C_p , \C_p) \; ,
\]
where $\Ch$ is the category $\Rep_{\pi_1 (X,x)} (\C_p)$, which makes the following diagram commutative:
\begin{equation}
  \label{eq:5}
\vcenter{  \xymatrix{
\Ext^1_{\eB_{X_{\C_p}}} (\Oh , \Oh) \ar[rr]^{\rho_*} \ar@{=}[d] & & \Ext^1_{\Ch} (\C_p , \C_p) \ar[d]^{\wr} \\
H^1 (X , \Oh) \otimes_{\oQ_p} \C_p \ar[rr]^{\rm Hodge-Tate} & & H^1_{\et} (X , \Q_p) \otimes_{\Q_p} \C_p
}}
\end{equation}
Here the right vertical isomorphism is defined as in the case of abelian varieties, and the lower horizontal map comes from the Hodge--Tate decomposition of $H^1_{\et} (X, \Q_p) \otimes \C_p$. 
\end{cor}

\begin{proof}
  According to \cite{De-We2}, Theorem 10, $\eB_{X_{\C_p}}$ is stable under extensions, so that every vector bundle on $X_{\C_p}$ which is an extension of the trivial bundle by itself lies in $\eB_{X_{\C_p}}$. By \cite{De-We2}, Proposition 21, the association $F \mapsto \rho_F$, mapping a vector bundle $F$ to the continuous representation $\rho_F : \pi_1 (X,x) \to \Aut_{\C_p} (F_x)$, respects exact sequences. Hence it induces a homomorphism on Yoneda $\Ext$ groups. Denote by $f : X \to A$ the morphism of $X$ into its Jacobian with $f (x) = 0$. By Lemma \ref{t6}, the middle square  in the following diagram is commutative:
\[
\xymatrix{
H^1 (A_{\C_p} , \Oh) \ar[r]^{\overset{f^*}{\sim}} \ar@{=}[d] & H^1 (X_{\C_p} , \Oh) \ar@{=}[d] \\
\Ext^1_{\eB_{A_{\C_p}}} (\Oh , \Oh) \ar[r]^{f^*} \ar[d]^{\rho_*} & \Ext^1_{\eB_{X_{\C_p}}} (\Oh , \Oh) \ar[d]^{\rho_*} \\
\Ext^1_{\Rep_{TA} (\C_p)} (\C_p , \C_p) \ar[r]^-{\tilde{f}_*} \ar[d]^{\wr} & \Ext^1_{\Ch} (\C_p , \C_p) \ar[d]^{\wr} \\
H^1_{\et} (A , \Q_p) \otimes \C_p \ar[r]^{\overset{f^*}{\sim}} & H^1_{\et} (X , \Q_p) \otimes \C_p
}
\]
The outer squares are also commutative. Since the Hodge--Tate decomposition is functorial, our claim follows from Theorem \ref{t6}. \qed
\end{proof}
\section{Line bundles on abelian varieties}
\label{sec:4}

Let $A$ be an abelian variety with good reduction over $\oQ_p$ which is defined over the finite extension $K$ of $\Q_p$ in $\oQ_p$, i.e. $A = A_K \otimes_K \oQ_p$. Then for the dual abelian variety $\hA$ we have $\hA = \hA_K \otimes_K \oQ_p$. Let $\Ah$ be an abelian scheme with generic fibre $A$ and $\hAh$ its dual abelian scheme.

By Theorem \ref{t4}, all line bundles algebraically equivalent to zero on $A_{\C_p}$ are contained in the category $\eB_{A_{\C_p}}$. For every line bundle $L$ on $A_{\C_p}$, we have $\Aut_{\C_p} (L_x) = \C_p^*$. Hence $\rho$ induces a map
\[
\alpha : \hA (\C_p) \longrightarrow \Hom_c (TA , \C^*_p) \; ,
\]
where $\Hom_c (TA , \C^*_p)$ denotes the continuous homomorphisms from $TA$ to $\C^*_p$. Note that $\Hom_c (TA , \C^*_p) = \Hom_c (TA , \eo^*)$ since $\C^*_p = p^{\Q} \times \eo^*$. By Theorem \ref{t4} and Proposition \ref{t5}, $\rho$ is compatible with tensor products and $G_K = \Gal (\oQ_p / K)$-action. Hence $\alpha$ is a $G_K$-equivariant homomorphism.

We endow the $\Ge_m$-torsor associated to the Poincar\'e bundle over $\Ah \times \hAh$ with the structure of a biextension of $\Ah$ and $\hAh$ by $\Ge_m$, so that we have $\hAh = \uExt^1 (\Ah , \Ge_m)$ in the flat topology. For all $N \ge 1$ we denote by $\Ah_N$ and $\hAh_N$ the subschemes of $N$-torsion points. The long exact $\Hom / \Ext$-sequence associated to the exact sequence
\[
0 \longrightarrow \Ah_N \longrightarrow \Ah \xrightarrow{N} \Ah \longrightarrow 0
\]
induces an isomorphism $\hAh_N \silo \uHom (\Ah_N , \Ge_m)$. Hence we get for all $n \ge 1$ a homomorphism
\[
\hAh_N (\eo_n) \longrightarrow \Hom (\Ah_N (\eo_n) , \eo^*_n) \; .
\]
By composition with the reduction map
\[
A_N (\oQ_p) \longrightarrow \Ah_N (\oZ_p) \longrightarrow \Ah_N (\eo_n)
\]
and with the projection $TA \to A_N (\oQ_p)$, we get a homomorphism
\[
\hAh_N (\eo_n) \longrightarrow \Hom_c (TA , \eo^*_n) \; .
\]
Note that here $\eo_n$ carries the discrete topology. For $N \tei M$ the corresponding maps are compatible with the inclusion $\hAh_N (\eo_n) \hookrightarrow \hAh_M (\eo_n)$. 

Note that the abelian group $\hAh (\eo_n)$ is torsion since it is the union of the finite groups $\hAh_{\eo_K} (\eo_L / p^n \eo_L)$, where $L$ runs over the finite extensions of $K$ in $\oQ_p$. Hence we get a homomorphism $\hAh (\eo_n) \to \Hom_c (TA , \eo^*_n)$. Composition with the reduction map $\hA (\C_p) = \hAh (\eo) \to \hAh (\eo_n)$ induces a homomorphism
\[
\alpha_n : \hA (\C_p) \longrightarrow \Hom_c (TA , \eo^*_n) \; .
\]

\begin{theorem}
  \label{t9}
For every $n \ge 1$ and all $\ha \in \hA (\C_p)$ the $\eo^*_n$-valued character $\alpha_n (\ha)$ is the reduction of the $\eo^*$-valued character $\alpha (\ha)$ modulo $p^n$. Hence the homomorphism $\alpha : \hA (\C_p) \to \Hom_c (TA , \eo^*)$ is the inverse limit of the $\alpha_n$.
\end{theorem}

\begin{proof}
  We denote by $L$ the line bundle on $\Ah_{\eo}$ corresponding to the point $\ha \in \hA (\C_p) = \hAh (\eo)$. If $N$ is big enough, the reduction $\ha_n$ of $\ha$ modulo $p^n$ lies in $\hAh_N (\eo_n)$. By $L_n$ we denote the reduction of $L$ modulo $p^n$, i.e. $L_n$ is the line bundle on $\Ah_n$ corresponding to $\ha_n$. Then $N^* L_n$ is trivial on $\Ah_n$. Since we identified $\hAh$ with $\uExt^1 (\Ah , \Ge_m)$, the $\Ge_m$-torsor $\tL_n = L_n \ohne \{ \mbox{zero section} \}$ on $\Ah_n$ is endowed with the structure of an extension of $\Ah_n$ by $\Ge_{m , \eo_n}$. Moreover, with this identification the inclusion 
\[
i : \uHom (\Ah_N , \Ge_m) \simeq \hAh_N \hookrightarrow \hAh
\]
is given by pushout with respect to the exact sequence $0 \to \Ah_N \to \Ah \xrightarrow{N} \Ah \to 0$. Denote by $\varphi : \Ah_{n,N} \to \Ge_{m , \eo_n}$ the homomorphism corresponding to $\ha_n \in \hAh_N (\eo_n)$. Then $\tL_n$ is given by the following pushout diagram:
\begin{equation}
  \label{eq:6}
  \vcenter{
\xymatrix{
0 \ar[r] & \Ah_{n,N} \ar[r] \ar[d]_{\varphi} & \Ah_n \ar[r]^N \ar[d]_s & \Ah_n \ar[r] \ar@{=}[d] & 0 \\
0 \ar[r] & \Ge_{m,\eo_n} \ar[r] & \tL_n \ar[r] & \Ah_n \ar[r] & 0
}}
\end{equation}
By definition, $\alpha_n (\ha)$ is the map
\[
\alpha_n (\ha) : TA \longrightarrow A_N (\oQ_p) = \Ah_N (\oZ_p) \longrightarrow \Ah_N (\eo_n) = \Ah_{n,N} (\eo_n) \xrightarrow{\varphi} \eo^*_n \; .
\]
On the other hand, the reduction of
\[
\alpha (\ha) : TA \longrightarrow A_N (\oQ_p) \longrightarrow \eo*
\]
modulo $p^n$ is obtained from the map
\[
A_N (\oQ_p) \longrightarrow \eo^*_n
\]
associating to $a \in A_N (\oQ_p) = \Ah_N (\oZ_p)$ the element in $\eo^*_n$ corresponding to the natural action of $a_n \in \Ah_{n,N} (\eo_n)$ on $\Gamma (\Ah_n , N^* L_n) \silo \eo_n$. Here we can as well regard the natural action of $a_n \in \Ah_{n,N} (\eo_n)$ on $\Gamma (\Ah_n , N^* \tL_n) \silo \eo^*_n$ in the setting of $\Ge_m$-torsors.

Now the homomorphism $s : \Ah_n \to \tL_n$ from diagram (\ref{eq:6}) induces an element
\[
s_0 = (s, \id) : \Ah_n \longrightarrow \tL_n \times_{\Ah_n , N} \Ah_n = N^* \tL_n \quad \mbox{in} \; \Gamma (\Ah_n , N^* \tL_n)
\]
which is mapped to $s_0 \verk t_{a_n}$ via the action of $a_n$, where $t_{a_n}$ denotes translation by $a_n$ on $\Ah_n$. By diagram (\ref{eq:6}), the corresponding element in $\eo^*_n$ is equal to $\varphi (a_n)$. Therefore $\alpha (\ha)$ reduces to $\alpha_n (\ha)$ modulo $p^n$.
\qed
\end{proof}

In \cite{Ta}, \S\,4, Tate considers the homomorphism
\[
\alpha_T : \hAh (p) (\eo) \longrightarrow \Hom_c (T (\Ah (p)) , U_1)
\]
defined by duality of the $p$-divisible groups $\Ah (p)$ and $\hAh (p)$. Here $U_1$ denotes the group of units congruent to $1$ in $\eo$. It follows from theorem \ref{t13} that $\alpha_T$ coincides with the restriction of $\alpha$ to the open subgroup $\hAh (p) (\eo)$ of $\hA (\C_p)$. 

We now consider the restriction $\alpha_{\tors}$ of the map $\alpha$ to the torsion part of $\hA (\C_p)$:
\[
  \alpha_{\tors} : \hA (\C_p)_{\tors} = \hA (\oQ_p)_{\tors} \longrightarrow \Hom_c (TA , \eo^*)_{\tors} = \Hom_c (TA , \mu) \; .
\]
Here $\mu = \mu (\oQ_p)$ is the group of roots of unity in $\eo^*$ or $\oQ^*_p$. Note that the Kummer sequence
on $A_{\et}$ induces an isomorphism
\[
i_A : H^1 (A , \mu_N) \silo H^1 (A , \Ge_m)_N \; .
\]

\begin{prop}
  \label{t10}
The map $\alpha_{\tors}$ is an isomorphism. On $\hA_N (\oQ_p)$ it coincides with the composition:
\[
\hA_N (\oQ_p) = H^1 (A , \Ge_m)_N \xrightarrow{i^{-1}_A \atop \sim} H^1 (A , \mu_N) = \Hom_c (TA , \mu_N) \; .
\]
\end{prop}

\begin{proof}
By Theorem \ref{t9}, the restriction of $\alpha$ to $\hA_N (\oQ_p)$ is the map
\[
\hA_N (\oQ_p) \longrightarrow \Hom (A_N (\oQ_p) , \mu_N ) = \Hom (TA , \mu_N)
\]
coming from Cartier duality $\hA_N \simeq \uHom (A_N , \Ge_m)$ over $\oQ_p$. The canonical identification $\Hom_c (TA , \mu_N) = H^1 (A , \mu_N)$ can be factorized by the isomorphisms
\[
\Hom (A_N , \mu_N) \silo \Ext^1 (A , \mu_N) \silo H^1 (A , \mu_N) \; ,
\]
where the first map is induced by the exact sequence $0 \to A_N \to A \to A \to 0$ and the second map is the forgetful map associating to an extension the corresponding $\mu_N$-torsor. Since the diagram
\[
\xymatrix{
\Hom (A_N , \mu_N) \ar[r]^{\sim} \ar[d]_{\wr} & \Ext^1 (A , \mu_N) \ar[r]^{\sim} \ar[d] & H^1 (A , \mu_N) \ar[d]^{i_A}_{\wr} \\
\Hom (A_N , \Ge_m) \ar[r]^{\sim} & \Ext^1 (A , \Ge_m)_N \ar[r]^{\sim} \ar@{=}[d] & H^1 (A , \Ge_m)_N \\
\hA_N (\oQ_p) \ar[u]_{\wr} \ar@{=}[r] & \hA_N (\oQ_p) & 
}
\]
commutes, our claim follows.  \qed
\end{proof}

Next, we need an elementary lemma. Consider an abelian topological group $\pi$ which fits into an exact sequence of topological groups
\[
0 \longrightarrow H \longrightarrow \pi \longrightarrow \hat{\Z}^n \longrightarrow 0
\]
where $H$ is a finite discrete group. Later $\pi$ will be the abelianized fundamental group of an algebraic variety. Applying the functor $\Hom_c (\pi , \, \underline{\ \ })$ to the exact sequence
\[
0 \longrightarrow \mu  \longrightarrow \eo^* \xrightarrow{\log} \C_p \longrightarrow 0
\]
we get the sequence
\begin{equation}
  \label{eq:7}
\small  0 \longrightarrow \Hom_c (\pi , \mu ) \longrightarrow \Hom_c (\pi , \eo^*) \xrightarrow{\log_*} \Hom_c (\pi , \C_p) \longrightarrow 0 \; .
\end{equation}

\begin{lemma}
  \label{t11}
a) The sequence (\ref{eq:7}) is exact.\\
b) In the topology of uniform convergence $\Hom_c (\pi , \eo^*)$ is a complete topological group. It contains $(\eo^*)^n$ as an open subgroup of finite index and hence acquires a natural structure as a Lie group over $\C_p$. Its Lie algebra is $\Hom_c (\pi , \C_p)$ and the logarithm map is given by $\log_*$.
\end{lemma}

\begin{proof}
a)  Since $\Hom_c (\pi , \underline{\ \ })$ is left exact, it suffices to show that $\log_*$ is surjective. As $\Hom_c (H , \C_p) = 0$ we only have to show surjectivity of $\log_*$ for $\pi = \hat{\Z}^n$, hence for $\pi = \hat{\Z}$.
We first prove that the injective evaluation map $\ev : \Hom_c (\hat{\Z} , \eo^*) \longrightarrow \eo^* \; , \; \ev (\varphi) = \varphi (1)$
is surjective. \\
Set $U_1 = \{ x \in \eo^* \tei |x-1| < 1 \}$ and $U_0 = \{ x \in \eo^* \tei |x-1| < p^{-\frac{1}{p-1}} \}$. The logarithm provides an isomorphism
\[
\log : U_0 \silo V_0 = \{ x \in \C_p \tei |x| < p^{-\frac{1}{p-1}} \} \; ,
\]
whose inverse is the exponential map. Therefore $U_0$ is a $\Z_p$-module and it follows that
$\ev : \Hom_c (\Z_p , U_0) \silo U_0$
is an isomorphism. We claim that \linebreak
$\ev : \Hom_c (\Z_p , U_1) \hookrightarrow U_1$ is an isomorphism as well.
Fix some $b$ in $U_1$. We construct a continuous map $\psi : \Z_p \to U_1$ with $\psi (1) = b$ as follows. There is some $N \ge 1$ such that $b^{p^N} \in U_0$. Hence there is a continuous homomorphism $\varphi : p^N \Z_p \longrightarrow U_1$ such that $\varphi (p^N \nu) = (b^{p^N})^{\nu}$ for all $\nu \in \Z$. Because of the decomposition $\eo^* = \mu_{(p)} \times U_1$ the group $U_1$ is divisible. Hence there is a homomorphism $\psi' : \Z_p \to U_1$ whose restriction to $p^N \Z_p$ equals $\varphi$. It follows that $\psi'$ is continuous as well. Because of $\psi' (1)^{p^N} = \psi' (p^N) = b^{p^N}$ there is a root of unity $\zeta \in \mu_{p^N}$ with $\psi' (1) = \zeta b$. Take the continuous homomorphism $\psi'' : \Z_p \to \mu_{p^{\infty}} \subset U_1$ with $\psi'' (1) = \zeta^{-1}$ and set $\psi = \psi' \cdot \psi''$.

The natural projection $\hat{\Z} \to \Z_p$ induces a commutative diagram
\[
\xymatrix{
\Hom_c (\Z_p , U_1) \ar@{^{(}->}[rr] \ar[dr]^{\sim}_{\ev} & & \Hom_c (\hat{\Z} , U_1) \ar@{^{(}->}[dl]^{\ev} \\
 & U_1 & .
}
\]
It follows that $\ev : \Hom_c (\hat{\Z} , U_1) \to U_1$ is an isomorphism as well. Using the decomposition $\eo^* = \mu_{(p)} \times U_1$ where $\mu_{(p)}$ carries the discrete topology we conclude that $\ev : \Hom_c (\hat{\Z} , \eo^*) \longrightarrow \eo^*$ is an isomorphism. Using the commutative diagram
\[
\xymatrix{
\Hom_c (\hat{\Z} , \eo^*) \ar[r]^{\log_*} \ar[d]^{\wr \, \ev} & \Hom_c (\hat{\Z} , \C_p) \ar@{=}[r] & \Hom_c (\Z_p , \C_p) \ar[d]^{\wr \, \ev} \\
\eo^* \ar@{>>}[rr]^{\log} & & \C_p
}
\]
we see that $\log_*$ is surjective for $\pi = \hat{\Z}$ and hence in general.\\
b) With the topology of uniform convergence, $\Hom_c (\pi , \eo^*)$ becomes a topological group. This topology comes from the inclusion of $\Hom_c (\pi , \eo^*)$ into the $p$-adic Banach space $C^0 (\pi , \C_p)$ of continuous functions from $\pi$ to $\C_p$ with the norm 
\[
\| f \| = \max_{\gamma \in \pi} |f (\gamma)| \; .
\]
Since $\Hom_c (\pi , \eo^*)$ is closed in $C^0 (\pi , \C_p)$ it becomes a complete metric space and hence it is a complete topological group. We now observe that the continuous evaluation map $\ev : \Hom_c (\hat{\Z} , \eo^*) \silo \eo^*$ is actually a homeomorphism. Let $x_n \to x$ be a convergent sequence in $\eo^*$ and let $\varphi_{\nu} , \varphi : \hat{\Z} \to \eo^*$ be the continuous homomorphisms with $\varphi_{\nu} (1) = x_{\nu}$ and $\varphi (1) = x$. Since $\Z^{\ge 1}$ is dense in $\hat{\Z}$ we get
\begin{eqnarray*}
  \| \varphi - \varphi_{\nu} \| & = & \max_{\gamma \in \hat{\Z}} |\varphi (\gamma) - \varphi_{\nu} (\gamma)| = \sup_{n \ge 1} |\varphi (n) - \varphi_{\nu} (n) | = \sup_{n \ge 1} |x^n - x^n_{\nu}| \\
& = & |x - x_{\nu}| \sup_{n \ge 1} \, |\sum^{n-1}_{i=0} x^i x^{n-i-1}_{\nu}| \le |x - x_{\nu}| \; .
\end{eqnarray*}
Hence $\varphi_{\nu}$ converges uniformly to $\varphi$. 

It follows that $\Hom_c (\hat{\Z}^n , \eo^*)$ and $(\eo^*)^n$ are isomorphic as topological groups. The exact sequence of topological groups
\[
0 \to \Hom_c (\hat{\Z}^n , \eo^*) \to \Hom_c (\pi , \eo^*) \to \Hom_c (H , \eo^*) = \Hom_c (H , \mu_{|H|})
\]
shows that $\Hom_c (\pi , \eo^*)$ contains $(\eo^*)^n$ as an open subgroup of finite index. Hence $\Hom_c (\pi , \eo^*)$ becomes a Lie group over $\C_p$. It is clear that the analytic structure depends only on $\pi$ and not on the choice of an exact sequence \\
$0 \to H \to \pi \to \hat{\Z}^n \to 0$ as above. The remaining assertions have to be checked for $\pi = \hat{\Z}^n$ and hence for $\pi = \hat{\Z}$ only where they are clear by the preceeding observations.\qed
\end{proof}

\begin{rem}
  The proof shows that the topologies of uniform and pointwise convergence on $\Hom_c (\pi , \eo^*)$ coincide.
\end{rem}

By \cite{Bou}, III, \S7, no.6 there is a logarithm map on an open subgroup $U$
of the $p$-adic Lie group $\hA (\C_p)$, mapping $U \rightarrow  \Lie \hA (\C_p)$, such that the kernel consists of the torsion points in $U$. Since $\hA (\C_p) / U$ is torsion (see Theorem 4.1 in \cite{Co2}), the logarithm has a unique extension to the whole Lie group $\hA (\C_p)$. It is surjective since the $\C_p$-vector
space $\Lie \hA (\C_p)$ is uniquely divisible. Therefore we have the exact sequence
\begin{equation}
  \label{eq:8}
  0 \longrightarrow \hA (\C_p)_{\tors} \longrightarrow \hA (\C_p) \xrightarrow{\log} \Lie \hA (\C_p) = H^1 (A , \Oh) \otimes \C_p \longrightarrow 0 \; .
\end{equation}
Using proposition \ref{t2} and lemma \ref{t11} we therefore get a commutative diagram with exact lines and $G_K$-equivariant maps
\begin{equation}
  \label{eq:9}
\def\objectstyle{\scriptstyle}
 \vcenter{\xymatrix@-1.2pc{
0 \ar[r] & \hA (\C_p)_{\tors} \ar[r] \ar[d]_{\alpha_{\tors}}^{\wr} & \hA (\C_p) \ar[r]^{\log} \ar[d]_{\alpha} & \Lie \hA (\C_p) \ar@{=}[r] \ar[d]^{\tilde{\alpha}} & H^1 (A , \Oh) \otimes_{\oQ_p} \C_p \ar[r] & 0 \\
0 \ar[r] & \Hom_c (TA , \mu) \ar[r] & \Hom_c (TA , \eo^*) \ar[r] & \Hom_c (TA , \C_p) \ar@{=}[r] & H^1_{\et} (A , \Q_p) \otimes \C_p \ar[r] & 0 
}}
\end{equation}
Here $\tilde{\alpha}$ is the map induced by $\alpha$ on the quotients. 

We will prove next that $\alpha$ is a $p$-adically analytic map of $p$-adic Lie groups. It follows that $\tilde{\alpha} = \Lie \alpha$. 

\begin{lemma}
  \label{t12}
Let $\beta : \hA (\C_p) \times TA \to \eo^*$ be the pairing induced by the homomorphism $\alpha : \hA (\C_p) \to \Hom_c (TA , \eo^*)$. Then $\beta$ is continuous. In particular, $\alpha$ is continuous.
\end{lemma}

\begin{proof}
  Denote by $r_n : \hAh (\eo) \to \hAh (\eo_n)$ the reduction map. Since the kernel of $r_n$ is $p$-adically open, it contains an open neighbourhood $W \subseteq \hA (\C_p)$ of zero.

Fix $(\ha, \gamma) \in \hA (\C_p) \times TA$ mapping to $z = \beta (\ha , \gamma)$. We show that the preimage of the open neighbourhood $z (1 + p^n \eo)$ is open. Let $N \ge 1$ be big enough so that $r_n (\ha)$ is contained in $\hAh_N (\eo_n)$. If $U$ denotes the kernel of the projection $TA \to A_N (\oQ_p)$, the neighbourhood $(\ha + W, \gamma + U )$ of $(\ha, \gamma )$ maps to $z (1 + p^n \eo)$ under $\beta$. Since the topology on $\Hom_c (TA, \eo^*)$ is the topology of pointwise convergence by the remark following lemma \ref{t11}, continuity of $\beta$ implies continuity of $\alpha$. \qed
\end{proof}

Let us fix some $\gamma \in TA$, and denote by $\psi_{\gamma}$ the induced homomorphism
\[
\psi_{\gamma} = \beta (-,\gamma) : \hA (\C_p) \longrightarrow \eo^* \; .
\]

\begin{prop}
  \label{t13}
$\psi_{\gamma}$ is an analytic map, hence a Lie group homomorphism.
\end{prop}

\begin{proof}
  We will briefly write $\psi = \psi_{\gamma}$ in this proof. It suffices to show that $\psi$ is analytic in a neighbourhood of the zero element $e_{\C_p} \in \hA (\C_p)$.

Let $e \in \hAh_{\eo} (\eo)$ be the zero section of $\hAh_{\eo}$. Since $\hAh_{\eo}$ is smooth over $\eo$, there is a Zariski open neighbourhood $U \subseteq \hAh_{\eo}$ of $e$ of the form
\[
U = \spec \eo [x_1 , \ldots , x_{m+r}] / (f_1 , \ldots , f_m)
\]
such that the matrix $\left( \frac{\partial f_i}{\partial x_{r+j}} (e) \right)_{i,j = 1 \ldots m}$ is invertible over $\eo$. 

By the theorem of implicit functions (see e.g. \cite{Col}, A.3.4), $U (\eo)$ contains an open neighbourhood $V$ of $e$ in the $p$-adic topology, such that the projection map $q : U (\eo) \subseteq \eo^{m+r} \to \eo^r$ given by $(x_1 , \ldots , x_{m+r}) \longmapsto (x_1 , \ldots , x_r)$ maps $e$ to $(0, \ldots , 0)$ and induces a homeomorphism
\[
q : V \longrightarrow V_1
\]
between $V$ and an open ball $V_1 \subseteq \eo^r$ around zero. This is an analytic chart around $e_{\C_p}$. A function $f$ on $V$ is locally analytic around $e_{\C_p}$, iff it induces a function on $V_1$ which coincides on a ball $V'_1 \subseteq V_1$ around $0$ with a power series in $x_1, \ldots , x_r$ converging pointwise on $V'_1$.

Since by \cite{Bou}, III, \S\,7, Prop. 10 and 11 the logarithm map on $\hA (\C_p)$ is locally around $e_{\C_p}$ an analytic isomorphism respecting the group structures, there exists an open subgroup $H$ of $\hA (\C_p)$ such that $(p^{\nu} H)_{\nu \ge 0}$ is a basis of open neighbourhoods of $e_{\C_p}$. By shrinking $V$ if necessary, we can assume that $V \subseteq H$.

For $n \ge 1$ we denote by $r_n$ as before the reduction map
\[
r_n : \hA (\C_p) = \hAh (\eo) \longrightarrow \hAh (\eo_n) = \hAh_n (\eo_n) \; ,
\]
where $\hAh_n = \hAh \otimes \eo_n$.

Since the kernel of $r_n$ is an open subgroup of $\hA (\C_p)$, it contains $p^{\nu_n} H$ for a suitable $\nu_n \ge 0$. Hence $p^{\nu_n} V$ is contained in the kernel of $r_n$, which implies that for all $x \in V$ the point $r_n (x)$ lies in the scheme $\hAh_{n , p^{\nu_n}}$ of $p^{\nu_n}$-torsion points in $\hAh_n$. 

The element $\gamma$ in $TA$ induces a point in $\Ah_{n, p^{\nu_n}} (\eo_n)$, whose image under the Cartier duality morphism
\[
\Ah_{n, p^{\nu_n}} (\eo_n) \longrightarrow \Hom (\hAh_{n , p^{\nu_n}} , \Ge_{m, \eo_n})
\]
we denote by $\psi_n$. Then $\psi (x)$ for $x \in V$ is by definition the element in $\eo^* \subseteq \eo$ satisfying $\psi (x) \equiv \psi_n (r_n (x)) \mod p^n$ for all $n$.

Let $U_n = U \otimes_{\eo} \eo_n = \spec \eo_n [x_1 , \ldots , x_{m+r}] / (\overline{f}_1 , \ldots , \overline{f}_m)$ be the reduction of the affine subscheme $U \subseteq \hAh_{\eo}$.
We write $\overline{f}$ for the reduction of a polynomial $f$ over $\eo$ modulo $p^n$.

Then $U_n \cap \hAh_{n , p^{\nu_n}} = \spec \eo_n [x_1 , \ldots , x_{m+r}] /  \ea$ for some ideal $\ea$ containing $(\of_1 , \ldots , \of_m)$. Since $\psi_n$ is an algebraic morphism, it is given on $U_n \cap \hAh_{n , p^{\nu_n}}$ by a unit in $\eo_n [x_1 , \ldots , x_{m+r}] / \ea$, which is induced by a polynomial $\overline{g}_n \in \eo_n [x_1 , \ldots , x_{m+r}]$. Let $g_n \in \eo [x_1 , \ldots , x_{m+r}]$ be a lift of $\overline{g}_n$. 

The implicit function theorem also implies that possibly after shrinking $V$ and $V_1$, we find power series $\theta_1 , \ldots , \theta_m \in \C_p [[ x_1 , \ldots , x_r]]$ converging in all points of $V_1 \subseteq \eo^r$, such that the map $V_1 \xrightarrow{q^{-1}} V \subseteq U (\eo) \subseteq \eo^{m+r}$ is given by
\[
(x_1 , \ldots , x_r) \longmapsto (x_1 , \ldots , x_r , \theta_1 (x_1 , \ldots , x_r) , \ldots , \theta_m (x_1, \ldots, x_r)) \; .
\]

For all $i = 1 , \ldots , m$ and all $n \ge 1$ let $h_{i,n} \in \C_p [x_1 , \ldots , x_r]$ be a polynomial satisfying $\theta_i (x) - h_{i,n} (x) \in p^n \eo$ for all $x \in V_1$. We can obtain $h_{i,n}$ by truncating $\theta_i$ suitably.

Then the map $
V_1 \silo V \xrightarrow{r_n} U_n (\eo_n) \cap \hAh_{n,p^{\nu_n}} (\eo_n) \xrightarrow{\psi_n} \eo^*_n$ maps the point
$
(x_1 , \ldots , x_r)$ to $ \overline{g}_n (\ox_1 , \ldots , \ox_r , \overline{\theta_1 (x_1 ,\ldots, x_r)} , \ldots , \overline{\theta_m (x_1, \ldots, x_r)}) $.

Hence for all $x = (x_1 , \ldots , x_r) \in V_1$ we have 
\[
\psi (q^{-1} (x)) - g_n (x_1 , \ldots , x_r , h_{1,n} (x_1, \ldots, x_r), \ldots , h_{m,n} (x_1, \ldots, x_r)) \in p^n \eo \; .
\]

Thus $\psi$ is the uniform limit of polynomials with respect to the coordinate chart $V_1$. This implies our claim. \qed
\end{proof}

\begin{cor}
  \label{t14}
The homomorphism $\alpha : \hA (\C_p) \longrightarrow \Hom_c (TA , \eo^*)$ is analytic.
\end{cor}

\begin{proof}
  Recall from the proof of lemma \ref{t11} that evaluation in $1 \in \hat{\Z}$ induces an isomorphism $\Hom_c (\hat{\Z} , \eo^*) \simeq \eo^*$. Hence $\Hom_c (TA , \eo^*) \simeq (\eo^*)^{2g}$ by evaluation in a $\hat{\Z}$-basis
 $\gamma_1 , \ldots , \gamma_{2g}$ of $TA$. This isomorphism induces the analytic structure on $\Hom_c (TA, \eo^*)$. Hence our claim follows from Proposition \ref{t13}. \qed
\end{proof}

Let us now determine the Lie algebra map induced by $\alpha$.

Recall from section \ref{sec:3} that the map
\[
\theta^*_A : H^1 (A_K , \Oh) \otimes_K \C_p \longrightarrow H^1_{\et} (A , \Q_p) \otimes \C_p
\]
coming from the Hodge--Tate decomposition of $H^1_{\et} (A , \Q_p) \otimes \C_p$ is the dual of a $\Z_p$-linear homomorphism
\[
\theta_A : T_p A \longrightarrow \omega_{\hAh} (\eo)
\]
defined using the universal vectorial extension of $\Ah$.

\begin{theorem}
  \label{t15}
We have $\Lie \alpha = \theta^*_A$.
\end{theorem}

\begin{proof}
  We give two proofs of this fact.\\
1) According to \cite{Ta}, \S\,4, the diagram
\[
\xymatrix{
\hAh (p) (\eo) \ar[r]^{\log} \ar[d]_{\alpha_T} & \Lie \hAh (p) \ar@{=}[r] \ar[d]^{\Lie \alpha_T} & H^1 (A , \Oh) \otimes_{\oQ_p} \C_p \\
\Hom_c (T (\Ah (p)) , U_1) \ar[r]^{\log_*} & \Hom_c (T (\Ah (p)) , \C_p) \ar@{=}[r] & H^1_{\et} (A , \Q_p) \otimes \C_p
}
\]
commutes. As we have seen, $\alpha_T$ is the restriction of $\alpha$ to $\hAh (p) (\eo)$. Combining Coleman's work in \cite{Co1} and \cite{Co3} \S\,4 with Fontaine's results, specifically \cite{Fo}, Proposition 11, it follows that $\Lie \alpha_T = \theta^*_A$. Hence $\Lie \alpha = \theta^*_A$.\\
2) The result can also be proved directly. Consider for $\gamma \in T_p (A)$ the analytic map $\psi_{\gamma} : \hA (\C_p) \longrightarrow \eo^*$ induced by $\alpha$. It induces a Lie algebra homomorphism
\[
\Lie \psi_{\gamma} : \Lie \hA (\C_p) \longrightarrow \Lie \eo^* \silo \eo \; ,
\]
where we identify $\Lie \eo^*$ with $\eo$ by means of the invariant differential $\frac{dT}{T}$ on $\Ge_{m, \eo} = \spec \eo [T, T^{-1}]$.

If suffices to show that for all $\gamma \in T_p A$ the map $\Lie \psi_{\gamma}$ is given by the invariant differential $\theta_A (\gamma)$, i.e. that
\[
\psi_{\gamma}^* \frac{dT}{T} = \theta_A (\gamma) \; .
\]
We use the notation from the proof of proposition \ref{t13}. Fix some $n \ge 1$ and recall the morphism
\[
\psi_n : U_n \cap \hAh_{n, p^{\nu_n}} = \spec \eo_n [x_1 , \ldots , x_{m+r}] / \ea \longrightarrow \Ge_{m, \eo_n}
\]
induced by a polynomial $\og_n \in \eo_n [x_1 , \ldots , x_{m+r}]$. This $\og_n$ defines a morphism
\[
\varphi^{(n)}_n : U_n \longrightarrow \Ge_{a , \eo_n} \; .
\]
Denote by $i : \Ge_{m, \eo} \rightarrow \Ge_{a, \eo}$ the obvious closed immersion, and by $i_n$ the induced morphism over $\eo_n$. The diagram 
\[
  \vcenter{\xymatrix{
U_n \ar[r]^{\varphi^{(n)}_n} & \Ge_{a , \eo_n} \\
U_n \cap \hAh_{n, p^{\nu_n}} \ar[r]^-{\psi_n} \ar@{^{(}->}[u]^{\kappa_n} & \Ge_{m, \eo_n} \ar@{^{(}->}[u]_{i_n}
}}
\]
is commutative, where $\kappa_n$ is  the canonical closed immersion. 

Similarly,  the lift $g_n$ of $\overline{g}_n$ to $\eo[x_1,\ldots, x_{m+r}]$ defines a morphism $\varphi^{(n)} : U \to \Ge_{a, \eo}$ over $\eo$, which induces an analytic map $\varphi^{(n)} : U (\eo) \longrightarrow \eo$.

The space of invariant differentials $\omega_{\hAh} (\eo)$ is an $\eo$-lattice in the $\C_p$-vector space $\omega_{\hA} (\C_p)$. For any analytic map $h : V \to \C_p$ we denote by $(dh)_e$ the  element in the cotangential space $\omega_{\hA} (\C_p) \simeq M_{e_{\C_p}} / M^2_{e_{\C_p}}$ given by the class of $h- h(e)$ modulo $M^2_{e_{\C_p}}$. Here $M_{e_{\C_p}}$ is the ideal of germs of analytic functions vanishing at $e_{\C_p}$. 

If $j : U \hookrightarrow \A^{m+r}_{\eo}$ denotes the obvious closed immersion in affine space, the fact that $\left( \frac{\partial f_i}{ \partial x_{j+r}} (e) \right)_{i,j = 1 \ldots m}$ is invertible implies that
\[
\omega_{\hAh} (\eo) = \Gamma (\spec \eo , e^* \Omega^1_{U / \eo})
\] 
is freely generated over $\eo$ by the differentials $j^* (dx_1)_e , \ldots , j^* (dx_r)_e$. 

For all $x \in V$ we have $i \verk \psi (x) \equiv i_n \verk \psi_n (r_n( x)) \mod p^n$, hence $i \verk \psi (x) - \varphi^{(n)} (x) \in p^n \eo \subseteq \eo$. Since  $\varphi^{(n)}$ is a polynomial map, we may assume by shrinking $V$, if necessary, that for {\em all} $n$ the function 
$(i \verk \psi - \varphi^{(n)}) \verk q^{-1}$ is given by a pointwise converging power series $\sum_{I = (i_1 , \ldots , i_r)} a^{(n)}_I x^{i_1}_1 \ldots x^{i_r}_r$ on the chart $V_1 \stackrel{q^{-1}}{\longrightarrow}V$, where $V_1 = (p^t \eo)^r$ for some $t \ge 0$.  For every multi-index $I$ this implies $p^{t (i_1 + \ldots + i_r)-n} a_I^{(n)} \in \eo$. 
Then 
\[d (i \verk \psi - \varphi^{(n)})_e = a^{(n)}_{(1 , 0 , \ldots , 0)} (j^* d{x_1})_e + \ldots + a^{(n)}_{(0 , \ldots , 0,1)} (j^* dx_r)_e \in p^{n-t} \omega_{\hAh} (\eo).\]

In particular for $n \geq t$ this implies $d (i \verk \psi)_e \in \omega_{\hAh} (\eo)$. Under the isomorphism \\
$\omega_{\hAh} (\eo) / p^n \omega_{\hAh} (\eo) \to \omega_{\hAh} (\eo_n)$ the element $(d \varphi^{(n)})_e \in \omega_{\hAh} (\eo)$ maps to $(d \varphi^{(n)}_n)_e \in \omega_{\hAh} (\eo_n)$.  Moreover the diagram above implies that
\[
\kappa^*_n (d \varphi^{(n)}_n)_e = (d (i_n \verk \psi_n))_e = \psi^*_n \left( \frac{dT}{T} \right)_e \; .
\]
Besides, we can assume that $\nu_n \ge n$, so that $\eo_n$ is annihilated by $p^{\nu_n}$. Then the exact sequence 
\[\omega_{\hAh_n} \stackrel{(p^{\nu_n})^*}{\longrightarrow} \omega_{\hAh_n} {\longrightarrow}  \omega_{\hAh_{n,p^{\nu_n}}} \longrightarrow 0\]
induces an isomorphism
$\kappa^*_n : \omega_{\hAh_n}  \stackrel{\sim}{\longrightarrow} \omega_{\hAh_{n,p^{\nu_n}}}$.
This implies 
\begin{eqnarray*}
  \psi^* \left( \frac{dT}{T} \right)_e = d (i \verk \psi)_e & \equiv & (d \varphi^{(n)})_e \mod p^{n-t} \omega_{\hAh} (\eo) \\
& \equiv & (\kappa^*_n)^{-1} \psi^*_n \left( \frac{dT}{T} \right)_e \mod p^{n-t} \omega_{\hAh} (\eo) \; .
\end{eqnarray*}
Now we take a closer look at the map $\theta_A$.

By definition, $\theta_A (\gamma) \equiv p^{\nu_n} b_{p^{\nu_n}} \mod p^n \omega_{\hAh} (\eo)$, where $b_{p^{\nu_n}} \in E (\eo)$ is an arbitrary lift of the $p^{\nu_n}$-torsion point $a_{p^{\nu_n}}$ in $\Ah (\eo)$ induced by $\gamma \in T_pA$ to the universal vectorial extension $E$.

Cartier duality $[\, , \,] : \Ah_{n , p^{\nu_n}} \times \hAh_{n, p^{\nu_n}} \to \Ge_{m,\eo_n}$ induces a homomorphism
\[
\tau_n : \Ah_{n,p^{\nu_n}} \longrightarrow \omega_{\hAh_{n,p^{\nu_n}}} \; ,
\]
given by $a \mapsto [a , -]^* \frac{dT}{T}$. Now we use an argument of Crew (see \cite{Cr}, section 1 and also \cite{Ch}, Lemma A.3) to show that $\theta_A (\gamma) \equiv (\kappa^*_n)^{-1} \tau_n (\overline{a_{p^{\nu_n}}}) \mod p^n \omega_{\hAh} (\eo)$, where $\overline{a_{p^{\nu_n}}} \in \Ah_n (\eo_n)$ is the reduction of $a_{p^{\nu_n}}$. 

Namely, by \cite{Ma-Me}, Chapter I, (2.6.2), the universal vectorial extension $E_n = E \otimes \eo_n$ of $\Ah_n$ is isomorphic to the pushout of the sequence 
\[0 \to \Ah_{n, p^{\nu_n}} \to \Ah_n \to \Ah_n \to 0\]
 by $(\kappa_n^*)^{-1} \verk \tau_n$. Hence we have a commutative diagram with exact lines
\[
\xymatrix{
0 \ar[r] & \Ah_{n , p^{\nu_n}} \ar[r] \ar[d]_{\tau_n} & \Ah_n \ar[dd]^f \ar[r]^{p^{\nu_n}} & \Ah_n \ar[r] \ar@{=}[dd] & 0 \\
 & \omega_{\hAh_{p^{\nu_n}}} & & & \\
0 \ar[r] & \omega_{\hAh_n} \ar[u]^{\kappa^*_n}_{\wr} \ar[r] & E_n \ar[r] & \Ah_n \ar[r] & 0 \; .
}
\]
Let $\overline{b_{p^{\nu_n}}}$ be the image of $b_{p^{\nu_n}}$ under the reduction map $E (\eo) \to E_n (\eo_n)$. Since multiplication by $p^{\nu_n}$ on $\Ah (\eo)$ is surjective, we can find some $\overline{c} \in \Ah_n (\eo_n)$ with $p^{\nu_n} \overline{c} = \overline{a_{p^{\nu_n}}}$. Then $f (\overline{c})$ differs from $\overline{b_{p^{\nu_n}}}$ by an element in $\omega_{\hAh_n} (\eo_n)$, which implies
  \begin{eqnarray*}
    p^{\nu_n} \overline{b_{p^{\nu_n}}} & = & f (p^{\nu_n} \overline{c}) \\
& = & (\kappa^*_n)^{-1} \tau_n (\overline{a_{p^{\nu_n}}}) \; ,
  \end{eqnarray*}
so that indeed 
\[
\theta_A (\gamma) \equiv (\kappa^*_n)^{-1} \tau_n (\overline{a_{p^{\nu_n}}}) \mod p^n \omega_{\hAh} (\eo) \; .
\]
By definition, $\tau_n (\overline{a_{p^{\nu_n}}}) = \psi^*_n \left( \frac{dT}{T} \right)_e$, so that for all $n$
{\begin{eqnarray*}
 \psi^* \left( \frac{dT}{T} \right)_e  & \equiv & (\kappa^*_n)^{-1} \psi^*_n \left( \frac{dT}{T} \right)_e \mod p^{n-t} \omega_{\hAh} (\eo) \\
& \equiv & \theta_A (\gamma) \mod p^{n-t} \omega_{\hAh} (\eo) \; ,
\end{eqnarray*}
which implies our claim.} \qed
\end{proof}

\begin{cor}
  \label{t16}
The map $\alpha$ is injective. 
\end{cor}

\begin{proof}
This follows from theorem \ref{t15} and diagram (\ref{eq:9}) since $\tilde{\alpha} = \Lie \alpha = \theta^*_A$ is injective.
\qed
\end{proof}

By Theorem \ref{t15}, the following diagram is commutative: 
\[
\def\objectstyle{\scriptstyle}
\def\labelstyle{\scriptstyle}
\xymatrix{
\hA (\C_p) \ar@{=}[r] & H^1 (A_{\C_p} , \Oh^*)^0 \ar[r]^{\log} \ar[d]_{\alpha} & H^1 (A_{\C_p} , \Oh) \ar@{=}[r] \ar[d]^{\theta^*_A} & \Lie \Pic^0_{A_K/K} (\C_p) \ar[d]^{\Lie \alpha} \\
 & \Hom_c (TA , \C^*_p) \ar[r]^{\log_*} & \Hom_c (TA , \C_p) \ar@{=}[r] & \Lie \Hom_c (TA , \C^*_p) \; .
}
\]
This is in a certain way analogous to the following diagram for the Lie group $U = \left\{ \left( 
    \begin{smallmatrix}
      1 & * \\ 0 & 1
    \end{smallmatrix} \right) \right\} \subset \GL_2$ which we derive from Theorem \ref{t7}:
\[
\def\objectstyle{\scriptstyle}
\def\labelstyle{\scriptstyle}
\xymatrix{
\Ext^1_{\eB_{A_{\C_p}}} (\Oh, \Oh) \ar@{=}[r] & H^1 (A_{\C_p} , U (\Oh)) \ar[r]^{\overset{\log = \id}{\sim}} \ar[d]_{\rho_*} & H^1 (A_{\C_p} , \Lie U (\Oh)) \ar@{=}[r] & H^1 (A_{\C_p} , \Oh) \ar[d]^{\theta^*_A} \\
 & \Hom_c (TA , U (\C_p)) \ar[r]^-{\log_* = \id} & \Hom_c (TA , \Lie U (\C_p)) \ar@{=}[r] & H^1_{\et} (A , \Q_p) \otimes \C_p \; .
}
\]
In the first diagram, the underlying group is $\Ge_m$, in the second it is $U \simeq \Ge_a$.

The next corollary was already observed by Tate in his context of $p$-divisible groups, \cite{Ta}, Theorem 3.

\begin{cor}
  \label{t17}
The map $\alpha$ induces an isomorphism of abelian groups
\[
\alpha : \hA_K (K) \silo \Hom_{c, G_K} (TA, \eo^*) \; .
\]
\end{cor}

\begin{proof}
  According to \cite{Ta}, Theorem 1 and 2 we have $H^0 (G_K , \C_p) = K$ and \\
$H^0 (G_K , \C_p (-1)) = 0$. Hence the Hodge--Tate decomposition and theorem \ref{t15} imply that $\tilde{\alpha} = \Lie \alpha$ induces an isomorphism:
\[
\tilde{\alpha} : H^1 (A_K , \Oh) \silo H^0 (G_K , H^1_{\et} (A , \Q_p) \otimes \C_p) = \Hom_{c, G_K} (TA , \C_p) \; .
\]
We have $H^0 (G_K , \hA (\C_p)) = \hA_K (K)$. This follows for example by embedding $\hA_K$ into some $\Pa^N$ over $K$ and using the corresponding result for $\Pa^N$. The latter is a consequence of the decomposition $\Pa^N = \A^N \amalg \ldots \amalg \A^0$ over $K$ and the equality $H^0(G_K , \C_p) = K$.
The corollary follows by applying the 5-lemma to the diagram of Galois cohomology sequences obtained from (\ref{eq:9}). \qed
\end{proof}

We next describe the image of the map $\alpha$ on $\hA (\C_p)$.

\begin{defn} \label{t18}
  A continuous character $\chi : TA \to \C^*_p$ is called smooth if its stabilizer in $G_K$ is open. The group of smooth characters of $TA$ is denoted by $Ch^{\infty} (TA)$.
\end{defn}

Note that we have
\[
Ch^{\infty} (TA) = \varinjlim_{L/K} \Hom_{c, G_L} (TA , \eo^*)
\]
where $L$ runs over the finite extensions of $K$ in $\oQ_p$. It is also the biggest $G_K$-invariant subset $S$ of $\Hom_c (TA , \C^*_p)$ such that the $G_K$-action on $S^{\delta}$ is continuous. Here $S^{\delta}$ is $S$ endowed with the discrete topology.

Replacing $A_K$ by $A_L$ in Corollary \ref{t17} we find that $\alpha$ induces an isomorphism
\[
\alpha : \hA (\oQ_p) \silo Ch^{\infty} (TA) \subset \Hom_c (TA , \eo^*) \; .
\]

Let $Ch (TA)$ be the closure of $Ch^{\infty} (TA)$ in $\Hom_c (TA , \eo^*)$ or equivalently in $C^0 (TA , \C_p)$. Then $Ch (TA)$ is also a complete topological group.

\begin{theorem} \label{t19}
  The map $\alpha$ induces an isomorphism of topological groups
\[
\alpha : \hA (\C_p) \silo Ch (TA) \; .
\]
\end{theorem}

\begin{proof}
  By lemma \ref{t12}, $\alpha$ is continuous. Hence
\[
\alpha (\hA (\C_p)) = \alpha (\overline{\hA (\oQ_p)}) \subset \overline{\alpha (\hA (\oQ_p))} = Ch (TA) \; .
\]
It now suffices to show that $\alpha$ is a closed map. Namely, because of
\[
Ch^{\infty} (TA) \subset \alpha (\hA (\C_p)) \subset Ch (TA)
\]
it will follow that $\alpha (\hA (\C_p)) = Ch (TA)$ and $\alpha$ will be a homeomorphism onto its image.

So let $Y \subset \hA (\C_p)$ be a closed set. Let $\alpha (y_n)$ for $y_n \in Y$ be a sequence which converges to some $\chi$ in $\Hom_c (TA , \eo^*)$. Since the map $\log_*$ in (\ref{eq:7}) is continuous in the uniform topologies it follows that $\Lie \alpha (\log y_n) = \log_* \alpha (y_n)$ converges to $\log_* \chi$. Because of the equality 
\[
\Hom_c (TA , \C_p) = \Hom_{\Z_p} (T_p A , \C_p)
\]
the topology of uniform convergence on this space coincides with its topology as a finite dimensional $\C_p$-vector space. The map $\Lie \alpha$ is a $\C_p$-linear injection by the Hodge--Tate decomposition and theorem \ref{t15}. Hence it is a closed injection and therefore the sequence $\log y_n$ converges. As $\log$ is a local homeomorphism there is a convergent sequence $\tilde{y}_n \in \hA (\C_p)$ with $\log \tilde{y}_n = \log y_n$. Writing $\tilde{y}_n = y_n + t_n$ with $t_n \in \hA (\C_p)_{\tors}$ we get $\alpha (\tilde{y}_n) = \alpha (y_n) + \alpha (t_n)$. The sequence $\alpha (y_n)$ converges by assumption and the sequence $\alpha (\tilde{y}_n)$ converges because $\alpha$ is continuous. Hence the sequence $\alpha (t_n) = \alpha_{\tors} (t_n)$ converges. The groups $\hA (\C_p)_{\tors}$ and $\Hom_c (TA , \mu )$ are the kernels of the locally topological homomorphisms $\log$ resp. $\log_*$. Hence they inherit the discrete topology from the $p$-adic topologies on $\hA (\C_p)$ resp. $\Hom_c (TA, \eo^*)$. Therefore the algebraic isomorphism $\alpha_{\tors}$ is trivially a homeomorphism and hence the sequence $t_n$ converges.  It follows that the sequence $y_n$ converges to some $y \in Y$. By continuity of $\alpha$ the sequence $\alpha (y_n)$ converges to $\alpha (y)$. Thus $\alpha (Y)$ is closed as was to be shown. \qed
\end{proof}

The following example was prompted by a question of Damian Roessler.

{\bf Example} Fix some $\sigma$ in $G_K$. Since $\alpha$ is $G_K$-equivariant we know that if $\ha \in \hA (\C_p)$ corresponds to the character $\chi : TA \to \eo^*$ then $\sigma (\ha)$ corresponds to $^{\sigma} \chi = \sigma \verk \chi \verk \sigma^{-1}_*$. Here $\sigma_*$ is the action on $TA$ induced by $\sigma$.

How about the character $\sigma \verk \chi : TA \to \eo^*$? Using theorem \ref{t19} we will now show that it also corresponds to an element of $\hA (\C_p)$ provided that $A_K$ has complex multiplication over $K$. For this, we have to check that in the $CM$ case, the subgroup $Ch (TA)$ is invariant under the homeomorphism $\chi \mapsto \sigma \verk \chi$ of $\Hom_c (TA , \eo^*)$. It suffices to show that $Ch^{\infty} (TA)$ is invariant. For $\chi$ in $Ch^{\infty} (TA)$, there is a finite normal extension $N /K $ such that $\chi$ is $G_N$-invariant, i.e. $\tau^{-1} \chi \tau_* = \chi$ for all $\tau$ in $G_N$. It follows that
\[
\tau^{-1} (\sigma \chi) \tau_* = \sigma (\sigma^{-1} \tau^{-1} \sigma \chi) \tau_* = \sigma \chi (\sigma^{-1} \tau^{-1} \sigma)_* \tau_* = \sigma \chi [\sigma , \tau]_*
\]
where we define the commutator by $[\sigma , \tau] = \sigma^{-1} \tau^{-1} \sigma \tau$. By the $CM$ assumption, the image of $G_K$ in the automorphism group of $TA$ is abelian. Hence $[\sigma , \tau]_*$ acts trivially on $TA$ and we have thus shown that $\tau^{-1} (\sigma \chi) \tau_* = \sigma \chi$ for all $\tau \in G_N$. Hence $\sigma \verk \chi$ lies in $Ch^{\infty} (TA)$. This proves the claim. 

For $\ha \in \hA (\C_p)$ let $\ha_{\sigma} \in \hA (\C_p)$ be the element corresponding to $\sigma \verk \chi$ via theorem \ref{t19}. By construction, the map $(\sigma , \ha) \mapsto \ha_{\sigma}$ determines a new action of $G_K$ on $\hA (\C_p)$. It seems to be a nice exercise in $CM$-theory to give an explicit description of this action. 
\section{Line bundles on varieties and their $p$-adic characters}
\label{sec:5}

Consider a smooth and proper variety $X_K$ over a finite extension $K$ of $\Q_p$. Varieties are supposed to be geometrically irreducible. We assume that $H^1 (X_K)$ has good reduction in the sense that the inertia group $I_K$ of $G_K$ acts trivially on the \'etale cohomology group $H^1 (X , \Q_l)$ for some prime $l \neq p$. Here $X = X_K \otimes \oK$.

The abelianization of the fundamental group $\pi_1 (X , x)$ is independent of the choice of a base point $x$ and will be denoted by $\pi^{\abb}_1 (X)$. It carries an action of the Galois group $G_K$ even if $X_K$ does not have a $K$-rational point. We will now attach a $p$-adic character of $\pi^{\abb}_1 (X)$ to any line bundle $L$ on $X_{\C_p}$ whose image in the N\'eron--Severi group of $X_{\C_p}$ is torsion. 

 It is known that $B_K := \Pic^0_{X_K / K}$ is an abelian variety over $K$. Its dual is the Albanese variety $A_K = \Abb_{X_K / K}$ of $X_K$ over $K$. We put $B = B_K \otimes_K \oQ_p$ and $A = A_K \otimes_K \oQ_p$.

Using the Kummer sequence and divisibility  of $\Pic^0_{X_K/K} (\oQ_p)$ one gets an exact sequence
\begin{equation}
  \label{eq:17}
  0 \longrightarrow B_N (\oQ_p) \longrightarrow H^1 (X , \mu_N) \longrightarrow NS (X)_N \longrightarrow 0
\end{equation}
for every $N \ge 1$. Since the N\'eron--Severi group of $X$ is finitely generated it follows that $T_l B = H^1 (X , \Z_l (1))$ for every $l$. Thus $B_K$ and hence also $A_K$ have good reduction. For sufficiently large $N$ in the sense of divisibility we have
\[
NS (X)_{\tors} = NS (X)_N \; .
\]
Applying $\Hom (\, \underline{\ \ } \, , \mu_N)$ to the exact sequence (\ref{eq:17}) and passing to projective limits therefore gives an exact sequence of $G_K$-modules:
\begin{equation}
  \label{eq:18}
  0 \longrightarrow \Hom (NS (X)_{\tors} , \mu) \longrightarrow \pi^{\abb}_1 (X) \longrightarrow TA \longrightarrow 0 \; .
\end{equation}
Here we have used the perfect Galois equivariant pairing coming from Cartier duality
\[
A_N (\oQ_p) \times B_N (\oQ_p) \longrightarrow \mu_N \; .
\]

For every prime number $l$ the pro-$l$ part of the sequence (\ref{eq:18}) splits continuously since $T_l A$ is a free $\Z_l$-module. Hence (\ref{eq:18}) splits continuously and applying $\Hom_c (\underline{\ \ } , \eo^*)$ we get an exact sequence of $G_K$-modules
\[
  0 \longrightarrow \Hom_c (TA , \eo^*) \longrightarrow \Hom_c (\pi^{\abb}_1 (X) , \eo^*) \longrightarrow NS (X)_{\tors} \longrightarrow 0 \; .
\]
We set
\[
\Hom^0_c (\pi^{\abb}_1 (X) , \eo^*) = \Ker (\Hom_c (\pi^{\abb}_1 (X) , \eo^*) \longrightarrow NS (X)_{\tors}) \; .
\]
Recall from section 3 the continuous injective homomorphism
\[
  \alpha : \Pic^0_{X_K/K} (\C_p) = \hA (\C_p) \hookrightarrow \Hom_c (TA , \eo^*) = \Hom^0_c (\pi^{\abb}_1 (X) , \eo^*) \; .
\]
Moreover $\alpha$ is a locally analytic homomorphism of $p$-adic Lie groups over $\C_p$.  Using theorem \ref{t15}, we see that $\alpha$ fits into a commutative diagram with exact lines:
\begin{equation}
  \label{eq:21}
\vcenter{ \def\objectstyle{\scriptstyle}
 \xymatrix@C=10pt{
0 \ar[r] & \Pic^0_{X_K/K} (\C_p)_{\tors} \ar[d]^{\wr}_{\alpha_{\tors}} \ar[r] & \Pic^0_{X_K/K} (\C_p) \ar@{^{(}->}[d]_{\alpha} \ar[r]^{\log} & \Lie \Pic^0_{X_K/K} (\C_p) \ar@{^{(}->}[d]^{\Lie \alpha} \ar@{=}[r] & H^1 (X_K, \Oh) \otimes_K \C_p \ar[r] & 0 \\
0 \ar[r] & \Hom^0_c (\pi^{\abb}_1 (X) , \mu ) \ar[r] & \Hom^0_c (\pi^{\abb}_1 (X) , \eo^*) \ar[r]^{\log_*} & \Hom_c (\pi^{\abb}_1 (X) , \C_p) \ar@{=}[r] & H^1 (X_{\et} , \Q_p) \otimes \C_p \ar[r] & 0 \; .
}}
\end{equation}

Here we have set
\begin{eqnarray*}
  \Hom^0_c (\pi^{\abb}_1 (X) , \mu) & = & \Hom^0_c (\pi^{\abb}_1 (X) , \eo^*) \cap \Hom_c (\pi^{\abb}_1 (X) , \mu ) \\
& = & \Hom^0_c (\pi^{\abb}_1 (X) , \eo^*)_{\tors} \; .
\end{eqnarray*}
Furthermore, note that:
\[
\Lie \Hom^0_c (\pi^{\abb}_1 (X) , \eo^*) = \Hom_c (\pi^{\abb}_1 (X) , \C_p) \; .
\]
The map $\Lie \alpha$ coincides with the inclusion map coming from the Hodge--Tate decomposition of $H^1 (X_{\et} , \Q_p) \otimes \C_p$. This follows from theorem \ref{t15} and the functoriality of this decomposition.

Set
\[
Ch^{\infty} (\pi^{\abb}_1 (X))^0 = \varinjlim_{L/K} \Hom^0_{c , G_L} (\pi^{\abb}_1 (X) , \eo^*)
\]
and let $Ch (\pi^{\abb}_1 (X))^0$ be its closure in $\Hom^0_c (\pi^{\abb}_1 (X) , \eo^*)$. We make similar definitions with the $^0$'s omitted.

It follows from theorem \ref{t19} that $\alpha$ induces a topological isomorphism of complete topological groups:
\[
  \alpha : \Pic^0_{X_K/K} (\C_p) \silo Ch (\pi^{\abb}_1 (X))^0 \; .
\]
We will now extend the domain of definition of $\alpha$ to $\Pic^{\tau}_{X_K/K} (\C_p)$. This is the group of line bundles on $X_{\C_p}$ whose image in $NS (X_{\C_p}) = NS (X)$ is torsion. We thus have an exact sequence 
\begin{equation}
  \label{eq:23}
  0 \longrightarrow \Pic^0_{X_K/K} (\C_p) \longrightarrow \Pic^{\tau}_{X_K/K} (\C_p) \longrightarrow NS (X)_{\tors} \longrightarrow 0 \; .
\end{equation}

\begin{theorem}
  \label{t20}
There is a $G_K$-equivariant map $\alpha^{\tau}$ which makes the following diagrams with exact lines commute:
\[
  \vcenter{\xymatrix@C=10pt{
0 \ar[r] & \Pic^0_{X_K/K} (\C_p) \ar[d]^{\alpha} \ar[r] & \Pic^{\tau}_{X_K/K} (\C_p) \ar[d]^{\alpha^{\tau}} \ar[r] & NS (X)_{\tors} \ar@{=}[d] \ar[r] & 0 \\
0 \ar[r] & \Hom^0_c (\pi^{\abb}_1 (X) , \eo^*) \ar[r] & \Hom_c (\pi^{\abb}_1 (X) , \eo^*) \ar[r] & NS (X)_{\tors} \ar[r] & 0 
}}
\]
and
\[
 \vcenter{\def\objectstyle{\scriptstyle}
 \xymatrix@C=10pt{
0 \ar[r] & \Pic^{\tau}_{X_K/K} (\C_p)_{\tors} \ar[d]^{\wr\,\alpha^{\tau}_{\tors}} \ar[r] & \Pic^{\tau}_{X_K/K} (\C_p) \ar[d]^{\alpha^{\tau}} \ar[r]^{\log} & \Lie \Pic^{\tau}_{X_K/K} (\C_p) \ar@{^{(}->}[d]^{\Lie \alpha^{\tau} = \Lie \alpha} \ar@{=}[r] & H^1 (X_K, \Oh) \otimes_K \C_p \ar[r] & 0 \\
0 \ar[r] & \Hom_c (\pi^{\abb} _1 (X), \mu ) \ar[r] & \Hom_c (\pi^{\abb}_1 (X) , \eo^*) \ar[r] & \Hom_c (\pi^{\abb}_1 (X) , \C_p) \ar@{=}[r] & H^1 (X_{\et} , \Q_p) \otimes \C_p \ar[r] & 0 \; .
}}
\]
The map $\alpha^{\tau}$ is an injective and locally analytic homomorphism of $p$-adic Lie groups. Its restriction $\alpha^{\tau}_{\tors}$ to torsion subgroups is the inverse of the Kummer isomorphism:
\[
\scriptstyle i_X : \Hom_c (\pi^{\abb}_1 (X) , \mu) = \varinjlim_{N} H^1 (X , \mu_N) \silo \varinjlim_N H^1 (X , \Ge_m)_N = \Pic^{\tau}_{X_K/K} (\C_p)_{\tors} \; .
\]
The map $\alpha^{\tau}$ induces a topological isomorphism of complete topological groups
\[
\alpha^{\tau} : \Pic^{\tau}_{X_K/K} (\C_p) \silo Ch (\pi^{\abb}_1 (X)) \; .
\]
\end{theorem}

\begin{proof}
  We first note that $\Pic_{X_K/K} (\oQ_p)_N = \Pic_{X_K/K} (\C_p)_N$ and $\Pic^{\tau}_{X_K/K} (\oQ_p)_N = \Pic^{\tau}_{X_K/K} (\C_p)_N$ because this holds for $\Pic^0$ and because $NS (X_{\C_p}) = NS (X)$.

As $\Pic^0_{X_K/K} (\C_p)$ is divisible, the sequence (\ref{eq:23}) gives a short exact sequence
\[
0 \longrightarrow \Pic^0_{X_K/K} (\C_p)_{\tors} \longrightarrow \Pic^{\tau}_{X_K/K} (\C_p)_{\tors} \longrightarrow NS (X)_{\tors} \longrightarrow 0 \; .
\]
We claim that the following diagram with exact lines commutes:
\begin{equation}
  \label{eq:26}
\def\objectstyle{\scriptstyle}
\def\labelstyle{\scriptstyle}
\vcenter{\xymatrix@C=10pt{
0 \ar[r] & \Pic^0_{X_K/K} (\C_p)_{\tors} \ar[d]^{\alpha_{\tors}} \ar[r] & \Pic^{\tau}_{X_K/K} (\C_p)_{\tors} \ar[d]^{i^{-1}_X} \ar[r] & NS (X)_{\tors} \ar@{=}[d] \ar[r] & 0 \\
0 \ar[r] & \Hom^0_c (\pi^{\abb}_1 (X) , \mu) \ar[r] & \Hom_c (\pi^{\abb}_1 (X) , \mu) \ar[r]  & NS (X)_{\tors} \ar[r] & 0 \; .
}}
\end{equation}
If we make the identifications explicit which define the maps of the left square, we see that on $N$-torsion it is the outer rectangle of the following diagram
\[
\def\objectstyle{\scriptstyle}
\def\labelstyle{\scriptstyle}
\xymatrix@C=5pt{
\Pic^0_{X_K / K} (\C_p)_N \ar@{=}[r] \ar[d]_{\alpha_{\tors}} \ar@{}[dr] |{\fbox{\tiny 1}}& \hA_N (\C_p) \ar[d]_{\alpha_{\tors}} \ar@{=}[r] \ar@{}[dr] |{\fbox{\tiny 2}} & B_N (\C_p) \ar@{^{(}->}[rr] & & H^1 (X , \Ge_m)_N \ar[d]^{\wr i^{-1}_X} \\
\Hom^0_c (\pi^{\abb}_1 (X) , \mu_N) \ar@{=}[r] & \Hom_c (TA , \mu_N) \ar@{=}[r] & \Hom (A_N , \mu_N) \ar@{=}[r] & B_N \subset H^1 (X , \Ge_m)_N \ar[r]^-{\overset{i^{-1}_X}{\sim}} & H^1 (X , \mu_N) \; .
}
\]
Now, \fbox{\tiny 1} is commutative by definition, and \fbox{\tiny 2} commutes since the restriction of $\alpha$ to $\hA_N (\C_p)$ is the map $\hA_N (\C_p) \to \Hom_c (TA , \mu_N) = \Hom (A_N , \mu_N)$ coming from Cartier duality. Hence the outer rectangle commutes as well.

The right square in diagram (\ref{eq:26}) is commutative since the second map in the exact sequence (\ref{eq:17}) is induced by $i_X$.

We now define $\alpha^{\tau}$ on \[
\Pic^{\tau}_{X_K/K} (\C_p) = \Pic^0_{X_K/K} (\C_p) + \Pic^{\tau}_{X_K/K} (\C_p)_{\tors}\] 
by setting it equal to $\alpha$ on $\Pic^0_{X_K/K} (\C_p)$ and to $i^{-1}_X$ on $\Pic^{\tau}_{X_K/K} (\C_p)_{\tors}$. This is well defined since by the commutativity of (\ref{eq:26}) the maps $\alpha$ and $i^{-1}_X$ agree on 
\[
\Pic^0_{X_K/K} (\C_p) \cap \Pic^{\tau}_{X_K/K} (\C_p)_{\tors} = \Pic^0_{X_K/K} (\C_p)_{\tors} \; . 
\]
The remaining assertions follow without difficulty. Note that $\Pic^{\tau}_{X_K/K} (\C_p)_{\tors}$ carries the discrete topology as a subspace of $\Pic^{\tau}_{X_K/K} (\C_p)$ since it is the kernel of the locally topological log map. \qed
\end{proof}
\newpage
\end{document}